\documentclass[a4paper,12pt,dvips]{article} 
\usepackage{graphicx}  
\usepackage{amssymb}
\usepackage{amscd} 
\usepackage{latexsym} 
\usepackage{array}
\usepackage{enumerate}

\usepackage{theorem}
\theoremstyle{plain}
\theorembodyfont{\itshape}
\newtheorem{thm}{Theorem}[section]
\newtheorem{prop}[thm]{Proposition}
\newtheorem{lem}[thm]{Lemma}

\theorembodyfont{\rmfamily}

\newtheorem{rem}[thm]{Remark}
\newtheorem{prob}[thm]{Problem}   
\makeatletter

\@addtoreset{equation}{section} 
\makeatother 
\def\dfrac#1#2{{\displaystyle\frac{#1}{#2}}} 
\def\Tr{\mathrm{Tr}}
\def\M{\mathcal{M}}
\def\K{\mathcal{K}}
\def\E{\mathcal{E}}
\def\R{\mathcal{R}}
\def\D{{\mit\Delta}} 
\def\PVI{\mathrm{P}_{\mathrm{VI}}} 

\def\C{\mathbb{C}}
\def\P{\mathbb{P}}
\def\Re{\mathrm{Re}}
\def\Im{\mathrm{Im}}
\def\RH{\mathrm{RH}} 
\def\-{\phantom{-}}
\def\ss{\scriptstyle} 
\begin{document}
\title{\Large\bf B\"acklund Transformations of the Sixth \\ 
Painlev\'e Equation in Terms of \\  
Riemann-Hilbert Correspondence}      
\author{Michi-aki Inaba\thanks{Faculty of Mathematics, Kyushu University, 
Hakozaki, Higashi-ku, Fukuoka 812-8581 Japan,  
E-mail address: inaba@math.kyushu-u.ac.jp},   
Katsunori Iwasaki\thanks{Faculty of Mathematics, Kyushu University, 
Hakozaki, Higashi-ku, Fukuoka 812-8581 Japan,  
E-mail address: iwasaki@math.kyushu-u.ac.jp} \, and  
Masa-Hiko Saito\thanks{Department of Mathematics, Faculty of Science, 
Kobe University, Rokko, Kobe 657-8501 Japan,  
E-mail address: mhsaito@math.kobe-u.ac.jp}} 
\date{July 7, 2003} 
\maketitle 
\begin{abstract} 
It is well known that the sixth Painlev\'e equation $\PVI$ admits 
a group of B\"acklund transformations which is isomorphic to the 
affine Weyl group of type $\mathrm{D}_4^{(1)}$. 
Although various aspects of this unexpectedly large symmetry have been 
discussed by 
many authors, there still remains a basic problem yet to be 
considered, that is, the problem of characterizing the B\"acklund 
transformations in terms of Riemann-Hilbert correspondence. 
In this direction, we show that the B\"acklund transformations 
are just the pull-back of very simple transformations on the 
moduli of monodromy representations by the Riemann-Hilbert  
correspondence. 
This result gives a natural and clear picture of the B\"acklund 
transformations.   
\end{abstract} 
%
\section{Introduction} \label{sec:intro} 
It is well known that the sixth Painlev\'e equation $\PVI$ 
admits a group of B\"acklund transformations which is 
isomorphic to the affine Weyl group of type $\mathrm{D}_4^{(1)}$. 
Various aspects of this unexpectedly large symmetry have 
been discussed by many authors, e.g. Okamoto \cite{Okamoto2}, 
Arinkin and Lysenko \cite{AR2}, Noumi and Yamada \cite{NY}. 
See also Conte and Musette \cite{CM}, Fokas and Yortsos \cite{FY}, 
Manin \cite{Manin}, Sakai \cite{Sakai}, Watanabe \cite{Watanabe} 
and others. 
However, there still seems to remain a basic problem yet to 
be considered with a special attention, namely, the problem of  
characterizing the B\"acklund transformations in terms of 
Riemann-Hilbert correspondence. 
This problem naturally arises from the work of 
Iwasaki \cite{Iwasaki3,Iwasaki4}, which exploited the standpoint 
of studying $\PVI$ based on a geometry of moduli spaces of monodromy 
representations via the Riemann-Hilbert correspondence; this    
standpoint had previously been hinted at in Iwasaki 
\cite{Iwasaki1,Iwasaki2}.  
\par   
The Riemann-Hilbert correspondence is a map from a moduli of 
Fuchsian differential equations to a moduli of monodromies, 
associating to each Fuchsian equation its monodromy representation.  
On the other hand, $\PVI$ is a differential equation on the moduli 
of Fuchsian equations that describes the isomonodromic deformation. 
Moreover, a B\"acklund transformation is a discrete transformation 
on the moduli of Fuchsian equations that commutes with the solution 
flow of $\PVI$. 
Therefore, it is very natural to ask what kind of discrete 
transformation on the moduli of monodromies is induced from a 
B\"acklund transformation by the Riemann-Hilbert correspondence. 
It is expected that the induced transformation looks simpler 
than the original, that is, a B\"acklund transformation 
looks simpler when viewed as a transformation on the moduli of 
monodromies. 
So the story should proceed in the other way round:  
Start with some simple transformations on the moduli of monodromies, 
pull them back via the Riemann-Hilbert correspondence to the moduli 
of Fuchsian equations and obtain the B\"acklund transformations. 
\par 
The aim of this paper is to verify that the speculation above is 
true: the B\"acklund transformations are just the pull-back of  
almost identity transformations on the moduli of monodromies, namely,  
of those transformations which may alter local monodromy data but  
do not change ``global monodromy data".   
The affine Weyl group structure of the B\"acklund transformations 
will also be transparent from our point of view. 
We hope that our result and viewpoint give a natural and clear 
picture of the B\"acklund transformations for $\PVI$. 
\par 
We remark that there is a geometric way to understand 
B\"acklund transformations of $\PVI$ by means of families of 
spaces of initial conditions constructed by Okamoto \cite{O1}.    
As explained in \cite{MMT,Shioda-Takano,Sakai,STT02}, 
there exists a family of open algebraic surfaces parametrized by 
the $4$-dimensional space $\K$ of local exponents, 
which correspond to the spaces of initial conditions of $\PVI$. 
Then the affine Weyl group $W(\mathrm{D}_4^{(1)})$ acts on $\K$ 
in a natural way and the actions can be lifted to
birational transformations of the total space of the family of 
Okamoto spaces. 
Saito and Umemura \cite{SU} pointed out that a B\"acklund 
transformation corresponding to a reflection of 
$W(\mathrm{D}_4^{(1)})$ is nothing but a flop whose center is
a family of $(-2)$-rational curves contained in Okamoto spaces 
lying over the reflection hyperplanes in $\K$. 
Moreover, Saito and Terajima \cite{STe} clarified the relation 
between $(-2)$-curves in Okamoto spaces and Riccati solutions of
$\PVI$. Since a flop of family of algebraic surfaces appears as a
simultaneous resolution of rational double points, one can expect 
that transformations of reflection type are related to the 
simultaneous resolutions of versal deformations of
rational double points. 
\par 
On the other hand, Arinkin and Lysenko \cite{AL1} introduced the 
moduli space of $SL(2)$-bundles with connections on $\P^1$ 
parametrized by the local exponents and describe the B\"acklund 
transformations in \cite{AR2}. 
>From the viewpoint of the theory of isomonodromic deformations 
of flat connections (Fuchsian connections), the moduli space 
should corresponds to Okamoto spaces.
Unfortunately, they treated the moduli space mostly as stack and
restricted the parameter space to the complement of all reflection 
hyperplanes in order to avoid the reducible connections.
In a forthcoming paper \cite{IIS},
we shall construct moduli spaces of stable parabolic connections
over $\P^1$ for all parameters.   
By using our moduli spaces, we can give a more geometric and 
conceptual picture of the Riemann-Hilbert correspondence and
the B\"acklund transformations of $\PVI$.  
\section{Sketch of the Main Result} \label{sec:main}
We roughly state the main result of this paper. 
Here we are content with a sketch of outline, since a complete 
statement is possible only after all the necessary ingredients 
are prepared. 
Detailed explanations will be supplied in the subsequent sections.   
\par 
The sixth Painlev\'e equation $\PVI$ is a second order nonlinear 
ordinary differential equation 
\begin{equation} \label{eqn:PVI}
\begin{array}{rcl}
q_{xx} &=& \dfrac{1}{2} 
\left(\dfrac{1}{q}+\dfrac{1}{q-1}+\dfrac{1}{q-x} \right) q_x^2 
- \left(\dfrac{1}{x}+\dfrac{1}{x-1}+\dfrac{1}{q-x} \right) q \\
& & \vspace{-0.2cm} \\
&+& \dfrac{q(q-1)(q-x)}{2x^2(x-1)^2} 
\left\{\kappa_4^2 - \kappa_1^2 \dfrac{x}{q^2} + 
\kappa_2^2 \dfrac{x-1}{(q-1)^2} + 
(1-\kappa_3^2) \dfrac{x(x-1)}{(q-x)^2} \right\} 
\end{array} 
\end{equation}
for an unknown function $q = q(x)$, with complex parameters 
$\kappa = (\kappa_1,\kappa_2,\kappa_3,\kappa_4)$. 
Equation $\PVI$ is written $\PVI(\kappa)$ when the dependence on 
parameters $\kappa$ should be indicated explicitly. 
In describing the B\"acklund transformations, it is convenient to 
think of the parameter space as an affine space 
\begin{equation} \label{eqn:parameter}
\K = \{\kappa = (\kappa_0,\kappa_1,\kappa_2,\kappa_3,\kappa_4) 
\in \C^5\,:\, 
2\kappa_0 + \kappa_1 + \kappa_2 + \kappa_3 + \kappa_4 = 1\}.  
\end{equation}  
\par
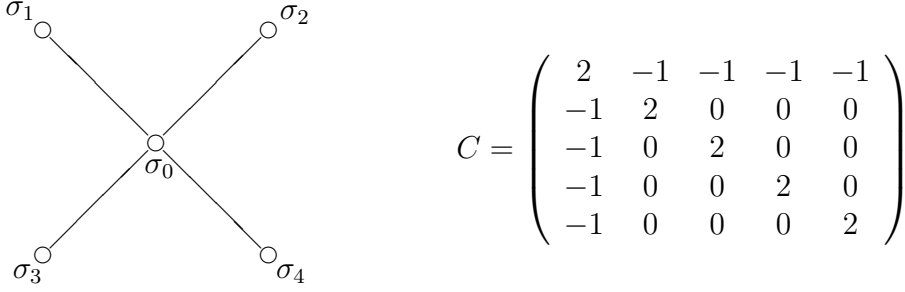
\begin{figure}[t] 
\begin{center}
\setlength{\unitlength}{1.0mm}
\begin{picture}(50,50)(40,5)
\put(10,10){\circle{2}}
\put(25,25){\circle{2}} 
\put(40,40){\circle{2}}
\put(40,10){\circle{2}}
\put(10,40){\circle{2}} 
\put(11,11){\line(1,1){13}}
\put(26,26){\line(1,1){13}}
\put(24,26){\line(-1,1){13}}
\put(26,24){\line(1,-1){13}}
\put(23.5,21){$\sigma_0$} 
\put(5,42){$\sigma_1$} 
\put(41.5,41.5){$\sigma_2$} 
\put(6,7){$\sigma_3$} 
\put(41,7){$\sigma_4$} 
\put(65,23){$C = \left(
\begin{array}{ccccc} 
\-2 &  -1 &  -1 &  -1 &  -1 \\
-1  & \-2 & \-0 & \-0 & \-0 \\
-1  & \-0 & \-2 & \-0 & \-0 \\
-1  & \-0 & \-0 & \-2 & \-0 \\
-1  & \-0 & \-0 & \-0 & \-2   
\end{array}
\right)$}    
\end{picture} 
\end{center}
\caption{Dynkin diagram and Cartan matrix of type $\mathrm{D}_4^{(1)}$} 
\label{fig:dynkin} 
\end{figure}

We realize the affine Weyl group of type $\mathrm{D}_4^{(1)}$ as 
an affine reflection group acting on $\K$, 
\[
W(\mathrm{D}_4^{(1)}) = \langle \sigma_0, \sigma_1, \sigma_2, \sigma_3, 
\sigma_4 \rangle, 
\]
where $\sigma_i$ is the reflection in the hyperplane $\kappa_i = 0$. 
In terms of the Cartan matrix $C = (c_{ij})$ of type $\mathrm{D}_4^{(1)}$ 
(see Figure \ref{fig:dynkin}), the reflection $\sigma_i$ is expressed as 
\begin{equation} \label{eqn:sigma} 
\sigma_i(\kappa_j) = \kappa_j - \kappa_i c_{ij}.    
\end{equation} 
A remarkable fact is that the affine Weyl group 
$W(\mathrm{D}_4^{(1)})$ lifts up to a transformation group of $\PVI$. 
Namely, each reflection $\sigma_i$ admits a lift $s_i$ that is a 
transformation of $\PVI(\kappa)$ to $\PVI(\sigma_i(\kappa))$, and 
the correspondence $\sigma_i \mapsto s_i$ induces an isomorphism 
between $W(\mathrm{D}_4^{(1)})$ and        
\begin{equation} \label{eqn:G}
G = \langle s_0, s_1, s_2, s_3, s_4 \rangle.   
\end{equation} 
The group $G$ is called the group of B\"acklund transformations 
for $\PVI$. 
The explicit form of $s_i$ will be given in \S\ref{sec:backlund} 
after some Hamiltonian formalisms for $\PVI$ are introduced in 
\S\ref{sec:hamilton}. 
\par 
We turn to a monodromy problem. 
Equation $\PVI(\kappa)$ is the isomonodromic deformation equation for 
a class of second order linear Fuchsian differential equations on 
$\P^1$ having four regular singular points with prescribed local 
exponents, where $\kappa$ is used to assign local exponents. 
We denote by $\E_t(\kappa)$ the moduli of relevant Fuchsian 
equations with regular singular points at $t = (t_1,t_2,t_3,t_4)$. 
The precise setting of $\E_t(\kappa)$ will be mentioned in 
\S \ref{sec:RH}.  
Let $\R_t(a)$ be the moduli of monodromy representations 
$\pi_1(\P^1\setminus \{t_1,t_2,t_3,t_4\}) \to SL_2(\C)$, up to 
Jordan equivalence, having prescribed local monodromy data  
$a = (a_1,a_2,a_3,a_4)$, where $a_i$ is defined to be the trace 
of the monodromy matrix $M_i$ around the singular point $t_i$, 
namely, 
\begin{equation} \label{eqn:ai} 
a_i = \mathrm{Tr}\, M_i \qquad (i = 1,2,3,4). 
\end{equation} 
Then the monodromy map or the Riemann-Hilbert correspondence  
\begin{equation} \label{eqn:RH1} 
\RH : \E_t(\kappa) \to \R_t(a)   
\end{equation} 
is defined by associating to each Fuchsian equation its monodromy 
representation class. 
In our setting, which will be detailed in \S\ref{sec:RH}, 
the correspondence of parameters $\kappa \mapsto a$ is given by 
\begin{equation} \label{eqn:ak} 
a_i = 2 \cos \pi \kappa_i \quad (i = 1,2,3), \qquad 
a_4 = -2 \cos \pi \kappa_4, 
\end{equation}
The minus sign for $a_4$ is not a misprint; $a_4$ is 
distinguished for a reason to be explained in \S\ref{sec:RH}.  
\par 
In Iwasaki \cite{Iwasaki3,Iwasaki4}, the representation space 
$\R_t(a)$ is realized as an affine cubic surface.  
Let us recall this construction. 
We introduce variables $x = (x_1,x_2,x_3)$ by     
\begin{equation} \label{eqn:xi}
x_i = \mathrm{Tr}(M_jM_k) \qquad \mathrm{for} \quad 
\{i,j,k\} = \{1,2,3\}.  
\end{equation} 
They are referred to as {\it global monodromy data}, since 
they carry a global information about monodromy; the product  
$M_jM_k$ is the monodromy matrix along a ``global'' loop 
surrounding the two singular points $t_j$ and $t_k$ 
simultaneously.   
Let $f(x,\theta)$ be a polynomial defined by 
\begin{equation} \label{eqn:f} 
f(x,\theta) = x_1 x_2 x_3 + x_1^2 + x_2^2 + x_3^2 
- \theta_1 x_1 - \theta_2 x_2 - \theta_3 x_3 + \theta_4,
\end{equation}  
where the coefficients 
$\theta = (\theta_1,\theta_2,\theta_3,\theta_4)$ are given by 
\begin{equation} \label{eqn:theta}
\theta_i = 
\left\{
\begin{array}{ll}
a_i a_4 + a_j a_k  & \qquad (i = 1,2,3), \\ 
\vspace{-0.2cm} & \\
a_1a_2a_3a_4 + a_1^2 + a_2^2 + a_3^2 + a_4^2 - 4  & \qquad  
(i = 4). 
\end{array} 
\right.
\end{equation}  
Then the representation space $\R_t(a)$ can be identified with 
an affine cubic surface 
\begin{equation} \label{eqn:S} 
\mathcal{S}(\theta) = \{\,x = (x_1,x_2,x_3) \in \C^3 \,:\, 
f(x,\theta) = 0\}.  
\end{equation}
Therefore, the Riemann-Hilbert correspondence (\ref{eqn:RH1}) 
is recast into 
\begin{equation} \label{eqn:RH2}
\RH : \E_t(\kappa) \to \mathcal{S}(\theta), 
\end{equation}
where the correspondence of parameters $\kappa \mapsto \theta$ 
is defined through (\ref{eqn:ak}) and (\ref{eqn:theta}). 
As a solution to the Riemann-Hilbert problem, the map 
(\ref{eqn:RH2}) is an analytic isomorphism onto a Zariski 
open subset of $\mathcal{S}(\theta)$.   
The following simple but fundamental observation is due to 
Terajima \cite{Terajima}.  
\begin{lem} \label{lem:invariant} 
Viewed as functions of $\kappa$, the coefficients $\theta$ 
are $W(\mathrm{D}_4^{(1)})$-invariants.  
\end{lem} 
\par
The proof is just by calculations. 
Note that the local monodromy data $a = (a_1,a_2,a_3,a_4)$  
are invariants of the reflections $\sigma_1, \sigma_2, 
\sigma_3, \sigma_4$, but not of $\sigma_0$, so that the 
unexpected switch $a \mapsto \theta$ has been necessary 
to obtain invariants of all five reflections $\sigma_i$. 
In this sense, the ``true'' local monodromy data might be 
attributed to $\theta$ rather than to $a$. 
Thanks to Lemma \ref{lem:invariant}, each B\"acklund 
transformation $s \in G$ induces an automorphism $r$ of 
$\mathcal{S}(\theta)$ such that the diagram  
\begin{equation} \label{eqn:CD}
\begin{CD}
\E_t(\kappa) @> s >> \E_t(\sigma(\kappa))  \\
      @V \RH VV     @VV \RH V              \\          
\mathcal{S}(\theta) @>> r > \mathcal{S}(\theta) 
\end{CD}  
\end{equation}
is commutative, where $\sigma \in W(\mathrm{D}_4^{(1)})$ is the 
transformation of parameters $\kappa$ underlying the B\"acklund 
transformation $s$. 
Now the following natural question occurs to us:  
\begin{prob} \label{prob:r} 
What is the transformation $r$ ? 
\end{prob} 
\par 
As one may expect naturally, this problem has a very simple 
solution:   
\begin{thm}[Main Theorem] \label{thm:main} 
The transformation $r$ in diagram $(\ref{eqn:CD})$ is just 
the identity$;$ that is to say, the B\"acklund transformations  
are those transformations which cover the identity 
transformation on the moduli of monodromies  through the 
Riemann-Hilbert correspondence.      
\end{thm} 
\par
The main result is far from trivial to the effect that 
the B\"acklund transformations, viewed as automorphisms of 
the cubic surface $\mathcal{S}(\theta)$, are distinguished from 
many other automorphisms to be the {\it identity} .  
In fact, $\mathcal{S}(\theta)$ admits a large number of 
automorphisms, even if the class is limited to algebraic 
automorphisms. 
In this respect, Iwasaki \cite{Iwasaki3,Iwasaki4} showed 
that the Riemann-Hilbert correspondence transforms  
the nonlinear monodromy of $\PVI$ into a modular group 
action on $\mathcal{S}(\theta)$ realized as a polynomial 
automorphism group on it. 
So $\mathcal{S}(\theta)$ admits at least infinitely 
many algebraic automorphisms labeled by the elements 
of a modular group. 
It may be said that the Riemann-Hilbert correspondence is 
such a map that sends a transcendental object on the phase 
space of $\PVI$ (the nonlinear monodromy of $\PVI$) to an 
algebraic object on the space of monodromies (a polynomial 
automorphism group on it), while collapsing an algebraic 
object on the former (the B\"acklund transformations) to 
a trivial object on the latter (the identity).  
\par 
In this paper we shall present an analytic proof of Theorem 
\ref{thm:main}, which is very simple in its essential idea 
(see \S\ref{sec:coals}) but requires some elaborate 
calculations in technicalities 
(see \S\ref{sec:accumulation}). 
It is desirable that there exists an alternative geometrical 
proof of the theorem. 
\par 
This paper is organized as follows. 
Some Hamiltonian formalisms for $\PVI$ are presented in 
\S\ref{sec:hamilton}, which enable us to describe the 
group of B\"acklund transformations in a symmetrical way 
in \S\ref{sec:backlund}.  
The Hamiltonian systems are characterized as isomonodromic 
deformation equations of second order Fuchsian differential 
equations with four regular singular points in 
\S\ref{sec:isomonod}.   
The B\"acklund transformations $s_1$, $s_2$, $s_3$ and $s_4$ 
are constructed as elementary gauge transformations in 
\S\ref{sec:gauge}. 
The Riemann-Hilbert correspondence is formulated in 
\S\ref{sec:RH}. 
With these preliminaries, our main result, Theorem 
\ref{thm:main}, is established in \S\ref{sec:coals}. 
The proof is based on the idea of coalescence of 
regular singular points along isomonodromic 
deformation. 
Here a key observation is that, by a coalescence procedure, 
a Fuchsian equation with four singular points degenerates 
into a Fuchsian equation with three singular points, and 
the difference of the two local exponents at the 
coalescent singular point is an invariant of the 
B\"acklund transformations; see Lemma \ref{lem:key}.     
In order for the idea to work, a certain technical 
lemma is needed concerning accumulation points of  
trajectories of the Hamiltonian system. 
This lemma is established 
in \S\ref{sec:accumulation} (Lemma \ref{lem:accumulation2}).  
The finial section, \S\ref{sec:heuristics}, is devoted 
to heuristics on finding B\"acklund transformations 
from our point of view.  
It provides a new way of discovering the hidden  
B\"acklund transformation $s_0$. 
\section{Hamiltonian Systems} \label{sec:hamilton}
In discussing B\"acklund transformations and monodromy problems 
related to $\PVI$, it is convenient to use a Hamiltonian 
system equivalent to the original single equation (\ref{eqn:PVI}). 
For this purpose, we shall employ the following three systems:  
\begin{enumerate} 
\item[(H1)] a Hamiltonian system with single time variable; see 
(\ref{eqn:HS});
\vspace{-0.2cm} 
\item[(H4)] a completely integrable Hamiltonian system with four 
time variables; see 
(\ref{eqn:HS4});
\vspace{-0.2cm}
\item[(H3)] a completely integrable Hamiltonian system with three 
time variables; see (\ref{eqn:HS3}).  
\end{enumerate}
We will flexibly use one or another of them depending upon 
contexts and purposes. 
\par
The first system (H1) is the most conventional one which often 
appears in the literature: 
\begin{equation} \label{eqn:HS}
\dfrac{\partial q}{\partial x} =  \dfrac{\partial h}{\partial p}, 
\qquad 
\dfrac{\partial p}{\partial x} = -\dfrac{\partial h}{\partial q},  
\end{equation}
where the Hamiltonian $h = h(q,p,x,\kappa)$ is given by 
\[
\begin{array}{rcl}
x(x-1) h &=&  
q(q-1)(q-x) p^2 - \{(\kappa_3-1) q(q-1) + \kappa_1 (q-1)(q-x) 
+ \kappa_2 q(q-x)\} p \\
\vspace{-0.3cm} & & \\
&&+ \kappa_0(\kappa_0+\kappa_4) (q-x).    
\end{array}
\]
Indeed, equation (\ref{eqn:PVI}) is recovered from system 
(\ref{eqn:HS}) by eliminating the variable $p$. 
\par  
A more symmetric description is feasible in terms of the 
second system (H4): 
\begin{equation} \label{eqn:HS4}
\dfrac{\partial q}{\partial t_i} =  \dfrac{\partial H_i}{\partial p}, 
\qquad 
\dfrac{\partial p}{\partial t_i} = -\dfrac{\partial H_i}{\partial q},  
\qquad (i = 1,2,3,4),  
\end{equation}
with four time variables $t = (t_1,t_2,t_3,t_4)$, 
the Hamiltonians $H_i = H_i(q,p,t,\kappa)$ being given by  
\begin{equation} \label{eqn:hamiltonian4}
\begin{array}{rcl}
(t_{ij}t_{ik}t_{il}) H_i &=& (q_iq_jq_kq_l) p^2 
-\{(\kappa_i-1) q_jq_kq_l + \kappa_j q_iq_kq_l + 
\kappa_k q_iq_jq_l + \kappa_l q_iq_jq_k\} p \\
\vspace{-0.3cm} & & \\
& &+ \kappa_0 q_i \{(\kappa_i-1) q_i + (\kappa_j+\kappa_0) q_j 
+ (\kappa_k+\kappa_0) q_k + (\kappa_l+\kappa_0) q_l \},   
\end{array} 
\end{equation} 
with $q_i = q-t_i$, $t_{ij} = t_i-t_j$ and 
$\{i,j,k,l\} = \{1,2,3,4\}$. 
We remark that (H1) is recovered from (H4) by a symplectic 
reduction. 
Indeed, we can observe that the diagonal action on 
$t = (t_1,t_2,t_3,t_4)$ of the M\"obius transformations 
lifts symplectically up to system (\ref{eqn:HS4}) and the 
associated symplectic reduction takes (\ref{eqn:HS4}) into 
(\ref{eqn:HS}) having the cross ratio  
\begin{equation} \label{eqn:cross} 
x = \dfrac{(t_1-t_3)(t_2-t_4)}{(t_1-t_2)(t_3-t_4)} 
\end{equation}  
as the only time variable. 
System (H3) is obtained from (H4) by letting $t_4$ tend to 
infinity. 
Indeed, 
\begin{equation} \label{eqn:limit1}
H_i(q,p,t,\kappa) \to h_i(q,p,t,\kappa) \quad (i = 1,2,3), 
\qquad H_4(q,p,t,\kappa) \to 0 
\quad \mathrm{as} \quad t_4 \to \infty, 
\end{equation}  
and system (\ref{eqn:HS4}) reduces to the Hamiltonian system   
\begin{equation} \label{eqn:HS3}
\dfrac{\partial q}{\partial t_i} =  \dfrac{\partial h_i}{\partial p}, 
\qquad 
\dfrac{\partial p}{\partial t_i} = -\dfrac{\partial h_i}{\partial q},  
\qquad (i = 1,2,3).   
\end{equation}
with three time variables $t = (t_1,t_2,t_3)$, 
where the Hamiltonians $h_i = h_i(q,p,t,\kappa)$ are given by 
\begin{equation} \label{eqn:hamiltonian3}
(t_{ij}t_{ik}) h_i = 
(q_iq_jq_k) p^2-\{(\kappa_i-1)q_jq_k+\kappa_jq_kq_i+\kappa_kq_iq_j\}p 
+ \kappa_0(\kappa_0+\kappa_4) q_i,   
\end{equation}
with $\{i,j,k\} = \{1,2,3\}$. 
The symplectic reduction mentioned above amounts to taking the 
normalization $t_1 = 0$, $t_2 = 1$, $t_3 = x$, $t_4 = \infty$. 
Then system (\ref{eqn:HS3}) with $i = 3$ yields (\ref{eqn:HS}). 
\section{Group of B\"acklund Transformations} 
\label{sec:backlund}
Now we can state (recall) the explicit form of the B\"acklund 
transformations. 
For a symmetric description, we shall represent it in terms of the 
Hamiltonian system with four time variables (\ref{eqn:HS4}). 
Necessary modifications for the systems (\ref{eqn:HS}) and 
(\ref{eqn:HS3}) are a routine work. 
The B\"acklund transformation $s_i$ corresponding to the 
reflection $\sigma_i$ is expressed as 
\begin{equation} \label{eqn:si} 
\begin{array}{rll}
s_0 &\left\{
\begin{array}{rcl}
s_0(t_j) &=& t_j, \\
\vspace{-0.3cm} & & \\
s_0(q) &=& q + \dfrac{\kappa_0}{p}, \\ 
\vspace{-0.3cm} & & \\
s_0(p) &=& p,
\end{array}
\right. & \qquad (i = 0),  \\
 &  & \\
s_i &\left\{
\begin{array}{rcl}
s_i(t_j) &=& t_j, \\ 
\vspace{-0.3cm} & & \\
s_i(q) &=& q, \\ 
\vspace{-0.3cm} & & \\
s_i(p) &=& p - \dfrac{\kappa_i}{q-t_i}, 
\end{array}
\right. & \qquad (i = 1,2,3,4). 
\end{array} 
\end{equation} 
\par
Noumi and Yamada \cite{NY} expressed (\ref{eqn:si}) in a unified 
manner by introducing variables  
\begin{equation} \label{eqn:invariant} 
q_i = \left\{
\begin{array}{ll}
p     & \qquad (i = 0), \\
\vspace{-0.3cm} & \\
q-t_i & \qquad (i = 1,2,3,4). 
\end{array} 
\right. 
\end{equation}
We remark that the fourth variable $q_4$ was degenerating 
in \cite{NY}, since they employed system (\ref{eqn:HS}) as 
their representation of $\PVI$, for which $t_4 = \infty$ and 
hence $q_4 = \infty$ was not observable.  
In any case, in terms of the variables $q_i$, 
formula (\ref{eqn:si}) together with (\ref{eqn:sigma}) is 
expressed as  
\begin{equation} \label{eqn:backlund} 
s_i(\kappa_j) = \kappa_j - \kappa_i c_{ij}, \qquad 
s_i(t_j) = t_j, \qquad 
s_i(q_j) = q_j + \dfrac{\kappa_i}{q_i} u_{ij},   
\end{equation} 
where $C = (c_{ij})$ is the Cartan matrix indicated in Figure 
\ref{fig:dynkin} and $u_{ij}$ is defined by 
\[
u_{ij} = \{q_i, q_j\} \qquad (i,j = 0,1,2,3,4).  
\] 
Here $\{f,g\}$ denotes the Poisson bracket, 
\[
\{f,g\} = \dfrac{\partial f}{\partial p} \dfrac{\partial g}{\partial q} 
- \dfrac{\partial f}{\partial q} \dfrac{\partial g}{\partial p}.   
\]
Explicitly, the matrix $U = (u_{ij})$ is given by  
\[
U = \left(
\begin{array}{ccccc} 
\- 0 & 1 & 1 & 1 & 1 \\
-1   & 0 & 0 & 0 & 0 \\
-1   & 0 & 0 & 0 & 0 \\
-1   & 0 & 0 & 0 & 0 \\
-1   & 0 & 0 & 0 & 0    
\end{array} 
\right).  
\]
\begin{rem} \label{rem:s4} 
When $t_4 = \infty$, formula (\ref{eqn:si}) with $i = 4$ 
should be interpreted as 
\[
s_4(t_j) = t_j, \qquad s_4(q) = q, \qquad s_4(p) = p. 
\]  
Namely, the lift $s_4$ of the reflection $\sigma_4$ acts on 
$(t,q,p)$ trivially. 
\end{rem} 
\par
The B\"acklund transformations $s_1$, $s_2$, $s_3$, $s_4$ 
are not so difficult to understand; as will be seen in 
\S\ref{sec:gauge}, they are nothing other than elementary 
gauge transformations of Fuchsian differential equations. 
Much more difficult and hence intriguing is the transformation 
$s_0$, whose existence is a sort of mystery. 
Hence the main body of this paper is devoted to understanding 
$s_0$. 
\par 
The symmetry of $\PVI$ with respect to the B\"acklund 
transformations is stated as follows. 
\begin{thm} \label{thm:backlund} 
Under the action of the B\"acklund transformation group $G$, the 
Hamiltonian system $(\ref{eqn:HS4})$ is invariant, that is, the 
action of $G$ commutes with the Hamiltonian vector fields 
\begin{equation} \label{eqn:vectorfields}
X_i = \dfrac{\partial}{\partial t_i} + 
\dfrac{\partial H_i}{\partial p} \dfrac{\partial}{\partial q} - 
\dfrac{\partial H_i}{\partial q} \dfrac{\partial}{\partial p}. 
\qquad (i = 1,2,3,4) 
\end{equation} 
\end{thm}  
{\it Proof}.  
The proof is by a direct check of condition (\ref{eqn:indep}) 
in the following lemma. \hfill $\Box$
\begin{lem} \label{lem:backlund}
The transformation $s_i$ commutes with the vector field $X_j$ if 
and only if 
\begin{equation} \label{eqn:indep}
s_i(H_j) - H_j + \delta_{ij} \dfrac{\kappa_i}{q_i} 
\quad \mbox{is independent of} \quad (q, p), 
\end{equation}
where $\delta_{ij}$ is Kronecker's delta symbol. 
\end{lem}
{\it Proof}. 
We prove the lemma for $i = 1,2,3,4$. 
The case $i = 0$ can be treated in a similar manner. 
Write $\bar{t}_j = s_i(t_j)$, $\bar{q} = s_i(q)$, $\bar{p} = s_i(p)$ 
and $\bar{H}_j = s_i(H_j)$. 
The transformation rule (\ref{eqn:si}) reads     
\[
t_j = \bar{t}_j, \qquad 
q = \bar{q}, \qquad 
p = \bar{p} + \dfrac{\kappa_i}{\bar{q}_i},     
\] 
where $\bar{q}_i = \bar{q} - \bar{t}_i$. 
Applying the chain rule of partial differentiations, we have  
\[
\dfrac{\partial}{\partial \bar{t}_j} = 
\dfrac{\partial}{\partial t_j} + \delta_{ij} \dfrac{\kappa_i}{q_i^2} 
\dfrac{\partial}{\partial p}, \qquad 
\dfrac{\partial}{\partial \bar{q}} = 
\dfrac{\partial}{\partial q} - \dfrac{\kappa_i}{q_i^2} 
\dfrac{\partial}{\partial p}, \qquad 
\dfrac{\partial}{\partial \bar{p}} = \dfrac{\partial}{\partial p},  
\]
and hence  
\[
\begin{array}{rcl}
s_i(X_j) 
&=& 
\dfrac{\partial}{\partial\bar{t}_j} + 
\dfrac{\partial\bar{H}_j}{\partial\bar{p}}\dfrac{\partial}{\partial\bar{q}}- 
\dfrac{\partial\bar{H}_j}{\partial\bar{q}} \dfrac{\partial}{\partial\bar{p}}\\
\vspace{-0.3cm} & & \\
&=& 
\left(\dfrac{\partial}{\partial t_j} + \delta_{ij} \dfrac{\kappa_i}{q_i^2} 
\dfrac{\partial}{\partial p}\right) 
+ \dfrac{\partial \bar{H}_j}{\partial p} 
\left(\dfrac{\partial}{\partial q} - \dfrac{\kappa_i}{q_i^2} 
\dfrac{\partial}{\partial p}\right)
- \left(\dfrac{\partial\bar{H_j}}{\partial q} - 
\dfrac{\kappa_i}{q_i^2} 
\dfrac{\partial\bar{H_j}}{\partial p}\right)\dfrac{\partial}{\partial p} \\
\vspace{-0.3cm} & & \\
&=& 
\dfrac{\partial}{\partial t_j} + 
\dfrac{\partial\bar{H}_j}{\partial p} \dfrac{\partial}{\partial q} - 
\left(\dfrac{\partial\bar{H}_j}{\partial q} - 
\delta_{ij} \dfrac{\kappa_i}{q_i^2} \right) \dfrac{\partial}{\partial p} 
\end{array}
\]
Therefore, we have $s_i(X_j) = X_j$ if and only if 
\[
\dfrac{\partial\bar{H}_j}{\partial p} = \dfrac{\partial H_j}{\partial p}, 
\qquad 
\dfrac{\partial\bar{H}_j}{\partial q} - 
\delta_{ij} \dfrac{\kappa_i}{q_i^2} = \dfrac{\partial H_j}{\partial q}. 
\]
The last condition is equivalent to the assertion (\ref{eqn:indep}). 
The proof is complete. \hfill $\Box$ 
\section{Isomonodromic Deformation} 
\label{sec:isomonod}
We discuss isomonodromy problems related to $\PVI$. 
Corresponding to the Hamiltonian systems with four time 
variables (\ref{eqn:HS4}) and with three time variables 
(\ref{eqn:HS3}), we set up two classes of Fuchsian 
differential equations, namely, (\ref{eqn:fuchs4}) and 
(\ref{eqn:fuchs3}) respectively. 
Then we consider their isomonodromic deformations. 
The Riemann-Hilbert correspondence for (\ref{eqn:fuchs3}) 
will be formulated in \S\ref{sec:RH}; a similar formulation 
for (\ref{eqn:fuchs4}) is possible, but omitted. 
\par  
For the Hamiltonian system (\ref{eqn:HS4}), we consider 
second order Fuchsian differential equations   
\begin{equation} \label{eqn:fuchs4} 
\dfrac{d^2f}{dz^2} - u_1(z) \dfrac{df}{dz} + u_2(z)f = 0
\end{equation}
on $\P^1$ with six singular points at 
$z = t_1, t_2, t_3, t_4, q, \infty$ such that 
\begin{enumerate}
\item[(A1)] $z = t_i$ is a regular singular point with exponents 
$0$ and $\kappa_i$$\mathrm{;}$
\vspace{-0.2cm}
\item[(A2)] $z = q$ is an apparent singular point with exponents 
$0$ and $2$$\mathrm{;}$ 
\vspace{-0.2cm}
\item[(A3)] $z = \infty$ is a removable singular point, that is, 
(\ref{eqn:fuchs4}) can be converted into a differential 
equation without singular point at $z = \infty$ by some 
transformation of the form $f = z^{-\kappa_0} g$.   
\end{enumerate} 
Here we recall the notion of an apparent singular point: 
A regular singular point $q$ of a second order Fuchsian differential 
equation is said to be {\it resonant} if the difference of the two 
exponents at $q$ is an integer. 
A resonant singular point $q$ falls into two cases; one is the generic 
case where a solution basis at $q$ involves the logarithmic function,  
and the other is the nongeneric case without logarithmic term in any 
solution basis. 
In the latter case, $q$ is called an {\it apparent} singular point. 
We remark that a removable singular point mentioned in (A3) is 
nothing but an apparent singular point such that the difference of 
exponents is one. 
We also remark that the number $\kappa_0$ in (A3) is uniquely 
determined, since equation (\ref{eqn:fuchs4}) must satisfy 
Fuchs' relation 
\begin{equation} \label{eqn:fuchsrel}
2\kappa_0+\kappa_1+\kappa_2+\kappa_3+\kappa_4 = 1. 
\end{equation}   
Note that the affine linear relation (\ref{eqn:fuchsrel}) is 
a source of the parameter space $\K$ in (\ref{eqn:parameter}). 
\par 
By conditions (A1), (A2), (A3), the coefficients 
of equation (\ref{eqn:fuchs4}) must be of the form    
\begin{equation} \label{eqn:coeff4} 
u_1(z) = \dfrac{1}{z-q}+\sum_{i=1}^4 \dfrac{\kappa_i-1}{z-t_i},  
\qquad 
u_2(z) = \dfrac{p}{z-q}-\sum_{i=1}^4 \dfrac{H_i}{z-t_i},     
\end{equation} 
where $(q,p,t,\kappa)$ are free parameters, while $H_i$ is a 
function of $(q,p,t,\kappa)$, whose explicit form can be sought 
out by Frobenius' method in the theory of Fuchsian differential 
equations; the result is that $H_i = H_i(q,p,t,\kappa)$ must be 
the Hamiltonians (\ref{eqn:hamiltonian4}) of the system 
(\ref{eqn:HS4}).  
\par 
It is well known that Painlev\'e-type equations arise from 
isomonodromic deformations of linear ordinary differential 
equations; see e.g. Fuchs \cite{Fuchs}, Schlesinger \cite{Schlesinger}, 
Jimbo, Miwa and Ueno \cite{JMU}, Fokas and Its \cite{FI}. 
In the original case of $\PVI$, we have the following theorem.  
\begin{thm} \label{thm:isomonod4}
The isomonodromic deformation of the Fuchsian differential 
equations $(\ref{eqn:fuchs4})$ is described by the completely 
integrable Hamiltonian system $(\ref{eqn:HS4})$.  
\end{thm} 
Writing $\PVI$ as a Hamiltonian system is due to 
Malmquist \cite{Malmquist}, Okamoto \cite{Okamoto1} and others. 
Iwasaki \cite{Iwasaki1,Iwasaki2} exploited an intrinsic 
(topological) reasoning for the Hamiltonian structure. 
\par  
Similarly, for the Hamiltonian system (\ref{eqn:HS3}), we consider 
Fuchsian differential equation   
\begin{equation} \label{eqn:fuchs3} 
\dfrac{d^2f}{dz^2} - v_1(z) \dfrac{df}{dz} + v_2(z)f = 0
\end{equation} 
on $\P^1$ with five singular points at 
$z = t_1, t_2, t_3, q, \infty$ such that   
\begin{enumerate}
\item[(B1)] $z = t_i$ is a regular singular point with exponents 
$0$ and $\kappa_i$;  
\vspace{-0.2cm}
\item[(B2)] $z = q$ is an apparent singular 
point with exponents $0$ and $2$;  
\vspace{-0.2cm}
\item[(B3)] $z = \infty$ is a regular singular point with 
exponents $\kappa_0$ and $\kappa_4 + \kappa_0$. 
\end{enumerate}
By conditions (B1), (B2), (B3), the coefficients of equation 
(\ref{eqn:fuchs3}) must be of the form     
\begin{equation} \label{eqn:coeff3}
v_1(z) = \dfrac{1}{z-q}+\sum_{i=1}^3 \dfrac{\kappa_i-1}{z-t_i}, 
\qquad  
v_2(z) = \dfrac{p}{z-q}-\sum_{i=1}^3 \dfrac{h_i}{z-t_i},  
\end{equation} 
where $(q,p,t,\kappa)$ are free parameters, while  
$h_i = h_i(q,p,t,\kappa)$ are the Hamiltonians 
(\ref{eqn:hamiltonian3}) of the system (\ref{eqn:HS3}).    
In view of (\ref{eqn:limit1}), formula (\ref{eqn:coeff3}) is 
obtained by taking the limit $t_4 \to 0$ in (\ref{eqn:coeff4}). 
The counterpart of Theorem \ref{thm:isomonod4} for  
equation (\ref{eqn:fuchs3}) is stated as follows.    
\begin{thm} \label{thm:isomonod3}
The isomonodromic deformation of the Fuchsian differential 
equation $(\ref{eqn:fuchs3})$ is described by the completely 
integrable Hamiltonian system $(\ref{eqn:HS3})$.  
\end{thm}
\section{Gauge Transformations} \label{sec:gauge} 
For $i = 1,2,3,4$, we can easily construct the B\"acklund 
transformation $s_i$ as an elementary gauge transformation 
of the Fuchsian equation (\ref{eqn:fuchs4}). 
Indeed, by the gauge transformation,  
\begin{equation} \label{eqn:gauge0}
f  = (z-t_i)^{\kappa_i} \bar{f}  \qquad (i = 1,2,3,4), 
\end{equation}
equation (\ref{eqn:fuchs4}) is transformed into another 
Fuchsian equation, say,  
\[
\dfrac{d^2\bar{f}}{dz^2} - \bar{u}_1(z) \dfrac{d\bar{f}}{dz} 
+ \bar{u}_2(z)\bar{f} = 0, 
\]
whose coefficients $\bar{u}_1$ and $\bar{u}_2$ must be of 
the same form as (\ref{eqn:coeff4}), say, 
\begin{equation} \label{eqn:coeff4bar} 
\bar{u}_1(z) = \dfrac{1}{z-\bar{q}}+\sum_{j=1}^4 
\dfrac{\bar{\kappa}_j-1}{z-\bar{t}_j},  
\qquad 
\bar{u}_2(z) = \dfrac{\bar{p}}{z-\bar{q}}-\sum_{j=1}^4 
\dfrac{\bar{H}_j}{z-\bar{t}_j},      
\end{equation} 
where $\bar{H}_j = H_j(\bar{q},\bar{p},\bar{t},\bar{\kappa})$.  
On the other hands, substituting (\ref{eqn:gauge0}) into 
(\ref{eqn:fuchs4}) implies that   
\begin{equation} \label{eqn:coeff4bar2} 
\bar{u}_1 = u_1 - \dfrac{2 \kappa_i}{z-t_i}, \qquad 
\bar{u}_2 = u_2 - \dfrac{\kappa_i}{z-t_i} u_1 
+ \dfrac{\kappa_i(\kappa_i-1)}{(z-t_i)^2}. 
\end{equation} 
Then substituting (\ref{eqn:coeff4}) into 
(\ref{eqn:coeff4bar2}) and comparing the result with 
(\ref{eqn:coeff4bar}), we have 
\[
\bar{\kappa}_j = \kappa_j - \kappa_i c_{ij}, 
\qquad \bar{t}_j = t_j, \qquad \bar{q} = q, \qquad 
\bar{p} = p - \dfrac{\kappa_i}{q-t_i}, 
\]
which is none other than the transformation $s_i$ in 
(\ref{eqn:si}) for $i = 1,2,3,4$. 
\section{Riemann-Hilbert Correspondence}
\label{sec:RH}
Let us formulate the Riemann-Hilbert correspondence for 
Fuchsian equations of the form (\ref{eqn:fuchs3}). 
It is convenient to define it in such a manner that monodromy 
matrices have values in the special linear group $SL_2(\C)$. 
To this end, by applying the gauge transformation 
\begin{equation} \label{eqn:gauge1}
f = \phi F \qquad \mathrm{with} \qquad 
\phi = (z-q) \prod_{i=1}^3 (z-t_i)^{\kappa_i/2},      
\end{equation}
equation (\ref{eqn:fuchs3}) should be normalized into a Fuchsian 
equation 
\begin{equation} \label{eqn:Fuchs3} 
\dfrac{d^2F}{dz^2} - V_1(z) \dfrac{dF}{dz} + V_2(z)F = 0 
\end{equation} 
with singular points at $z = t_1$, $t_2$, $t_3$, $q$, $\infty$  
satisfying the following conditions:    
\begin{enumerate} 
\item[(C1)] $z = t_i$ is a regular singular point with exponents 
$\pm\kappa_i/2$;  
\vspace{-0.2cm}
\item[(C2)] $z = q$ is an apparent singular 
point with exponents $\pm1$;  
\vspace{-0.2cm}
\item[(C3)] $z = \infty$ is a regular singular point with 
exponents $(3\pm\kappa_4)/2$.  
\end{enumerate}
Note that the coefficients of equation (\ref{eqn:Fuchs3}) and 
those of (\ref{eqn:fuchs3}) are related by  
\[ 
V_1 = v_1 - 2 \dfrac{\phi'}{\phi}, \qquad 
V_2 = v_2 - \dfrac{\phi'}{\phi} v_1 + \dfrac{\phi''}{\phi}. 
\] 
Equation (\ref{eqn:Fuchs3}) is called the {\it normal form} 
of (\ref{eqn:fuchs3}) and the two are often identified. 
Then the (normalized) monodromy of (\ref{eqn:fuchs3}) will be 
defined to be the monodromy of its normal form 
(\ref{eqn:Fuchs3}). 
\par 
Let $\E_t(\kappa)$ be the set of all Fuchsian equations of the 
form (\ref{eqn:fuchs3}), or equivalently of their normal form 
(\ref{eqn:Fuchs3}), with a prescribed value of parameters 
$\kappa$ and location of singular points $t = (t_1,t_2,t_3,t_4)$, 
where $t_4 = \infty$. 
Associating $(q,p)$ to equation (\ref{eqn:fuchs3}) yields an 
identification 
\begin{equation} \label{eqn:chart} 
\E_t(\kappa) \cong (\C \setminus \{t_1,t_2,t_3\}) \times \C 
\,\,:\,\, (q,p)\mbox{-space},  
\end{equation} 
through which $\E_t(\kappa)$ is thought of as a complex manifold. 
In an algebro-geometrical framework developed by the authors \cite{IIS}, 
the space $\E_t(\kappa)$ is a Zariski open chart of a certain moduli 
space $\M_t(\kappa)$ of rank two stable parabolic bundles with 
Fuchsian connections on $\P^1$. 
\par 
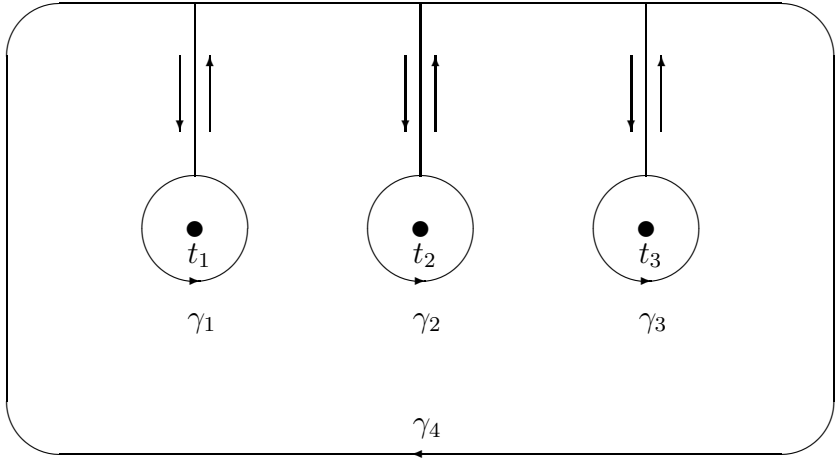
\begin{figure}[t] 
\begin{center}
\setlength{\unitlength}{1.0mm}
\begin{picture}(110,60)(-5,0)
\put(50,30){\oval(110,60)} 
\put(20,30){\circle*{2}}
\put(50,30){\circle*{2}} 
\put(80,30){\circle*{2}}
\put(20,30){\circle{14}}
\put(50,30){\circle{14}} 
\put(80,30){\circle{14}}
\put(20.5,23){\vector(1,0){0.2}} 
\put(50.5,23){\vector(1,0){0.2}} 
\put(80.5,23){\vector(1,0){0.2}} 
\put(49,0){\vector(-1,0){0.2}} 
\put(20,37){\line(0,1){23}}
\put(50,37){\line(0,1){23}}
\put(80,37){\line(0,1){23}}
\put(18,53){\vector(0,-1){10}}
\put(22,43){\vector(0,1){10}}
\put(48,53){\vector(0,-1){10}}
\put(52,43){\vector(0,1){10}}
\put(78,53){\vector(0,-1){10}}
\put(82,43){\vector(0,1){10}}
\put(19,25.5){$t_1$} 
\put(49,25.5){$t_2$}
\put(79,25.5){$t_3$}   
\put(19,17){$\gamma_1$} 
\put(49,17){$\gamma_2$} 
\put(79,17){$\gamma_3$}
\put(49,3){$\gamma_4$}  
\end{picture} 
\end{center}
\caption{The loops $\gamma_1$, $\gamma_2$, $\gamma_3$, $\gamma_4$}  
\label{fig:loops} 
\end{figure} 
We proceed to moduli of monodromies. 
Let $\gamma_i$ be a simple loop encircling the singular 
point $t_i$ as in Figure \ref{fig:loops}, with $\gamma_4$ being 
a loop around $t_4 = \infty$. 
Let $M_i$ be the monodromy matrix along the loop $\gamma_i$ 
of the Fuchsian equation (\ref{eqn:Fuchs3}).  
In view of (C1) and (C3), the matrix $M_i$ has the eigenvalues
\begin{equation}  \label{eqn:ev} 
\exp(\pm\pi\sqrt{-1}\kappa_i) \quad \mathrm{for} \quad i = 1,2,3; 
\qquad -\exp(\pm\pi\sqrt{-1}\kappa_4) \quad \mathrm{for} \quad 
i = 4. 
\end{equation} 
Hence $\det M_i = 1$, namely, $M_i \in SL_2(\C)$, and the 
trace $a_i = \mathrm{Tr}\, M_i$ is expressed as (\ref{eqn:ak}). 
Let $\R_t(a)$ be the space of monodromy representations 
of $\pi_1(\P^1 \setminus \{t_1,t_2,t_3,t_4 \})$ into 
$SL_2(\C)$, up to Jordan equivalence, whose monodromy matrices 
along the loop $\gamma_i$ have trace $a_i$. 
Then the Riemann-Hilbert correspondence 
$\RH : \E_t(\kappa) \to \R_t(a)$ in (\ref{eqn:RH1})  
is defined by associating to each Fuchsian equation 
(\ref{eqn:fuchs3}) the $SL_2(\C)$-monodromy representation 
class of its normal form (\ref{eqn:Fuchs3}),  
where two parameters $\kappa$ and $a$ are related by 
(\ref{eqn:ak}). 
Then, composed with a natural map 
$\R_t(a) \to \mathcal{S}(\theta)$, where $\mathcal{S}(\theta)$  
is the cubic surface in (\ref{eqn:S}),  
the Riemann-Hilbert correspondence is realized more concretely 
as a map into the cubic surface, namely, as the map 
$\RH : \E_t(\kappa) \to \mathcal{S}(\theta)$ in (\ref{eqn:RH2}).   
\par 
The construction above can be made relatively over the 
parameter spaces.   
Let $\E_t$ be the space of all Fuchsian equations of the 
form (\ref{eqn:fuchs3}), or equivalently of their normal 
form (\ref{eqn:Fuchs3}), where the regular singular points 
are fixed at $t = (t_1,t_2,t_3,t_4)$, $t_4 = \infty$, but 
the parameters $\kappa$ may vary.   
Let $\pi_1 : \E_t \rightarrow \K$ be the natural projection
associating to each equation in $\E_t$ its parameter $\kappa$.
Note that $\E_t(\kappa)$ is a fiber of this projection.
The space $\E_t$ is a Zariski open subset of the family
$\M_t$ of moduli spaces $\M_t(\kappa)$ over $\K$,
which is constructed in Inaba, Iwasaki and Saito \cite{IIS}.
Let $\R_t$ be the space of all monodromy representations
$\pi_1(\P^1\setminus \{t_1,t_2,t_3,t_4\}) \to SL_2(\C)$ up to
Jordan equivalence and let $\pi_2 : \R_t \to \C_a^4$ be the natural
projection associating to each representation in $\R_t$ its
local monodromy data $a$.  Recall that
$\R_t $ is the categorical quotient of the diagonal adjoint action
of $SL_2(\C)$ on $SL_2(\C)^{3}$ and that $\pi_2 : \R_t \to \C_a^4$ 
is a family of affine cubic surfaces defined by the equation  
$f(x, \theta(a)) = 0$ in  (\ref{eqn:f}) with the relation 
$\theta = \theta(a)$ in (\ref{eqn:theta}).  
Moreover, let
\[
\mathcal{S} = \{\,(x,\theta) \in \C^3 \times \C^4\,:\,
f(x,\theta) = 0\,\}, 
\]
and let $\pi_3 : \mathcal{S} \to \C_{\theta}^4$ be the natural
projection down to parameters $\theta$. 
Note that $ \pi:\R_t \to \C_a^4$ is obtained by the pullback 
of the family  $\pi_3 : \mathcal{S} \to \C_{\theta}^4$ by the 
finite morphism $\C_{a}^4 \to \C_{\theta}^4$.   
Now the (relative) Riemann-Hilbert correspondence $\RH$ is 
formulated as a commutative diagram:  
\[ 
\begin{CD} 
\M_t @> \RH >> \R_t @> \varphi >> \mathcal{S}  \\
 @V \pi_1 VV   @VV \pi_2 V  @VV \pi_3 V  \\
 \K  @>> u >      \C_{a}^4     @>> v > \C_{\theta}^4,  
\end{CD} 
\] 
where the maps $u : \kappa \mapsto a$ and 
$v : a \mapsto \theta$ are defined by (\ref{eqn:ak}) and 
(\ref{eqn:theta}); or as its contraction:   
\begin{equation} \label{eqn:RH3} 
\begin{CD} 
\M_t @> \RH >> \mathcal{S} \\
 @V \pi_1 VV   @VV \pi_3 V \\
 \K @>> vu > \C_{\theta}^4,  
\end{CD} 
\end{equation}  
\section{Coalescence of Regular Singular Points} 
\label{sec:coals}
A principal idea for establishing our main result, 
Theorem \ref{thm:main}, is to consider a coalescence of 
two regular singular points of Fuchsian equation, 
along an isomonodromic deformation.   
We first discuss the coalescence process only and then 
take the isomonodromic deformation into account. 
For this purpose we shall work with Fuchsian equation 
(\ref{eqn:fuchs3}). 
Of course, working with Fuchsian equation 
(\ref{eqn:fuchs4}) would lead us to the same conclusion; 
the latter choice would allow a more symmetrical discussion, 
but require somewhat heavier calculations. 
Here we employ (\ref{eqn:fuchs3}), preferring simpler 
calculations at the cost of a minor symmetry breaking. 
\begin{lem} \label{lem:coales} 
By the coalescence $t_k \to t_j$ of the singular points 
$t_j$ and $t_k$, equation $(\ref{eqn:fuchs3})$ with 
$(\ref{eqn:coeff3})$ degenerates into a Fuchsian equation 
with singular points at $z = t_i$, $t_j$, $q$, $\infty$,        
\begin{equation} \label{eqn:coalfuchs} 
\dfrac{d^2f}{dz^2} - w_1(z) \dfrac{df}{dz} + w_2(z)f = 0, 
\end{equation} 
whose coefficients $w_1(z)$ and $w_2(z)$ are expressed as   
\begin{eqnarray}
w_1(z) &=& \dfrac{1}{z-q} + \dfrac{\kappa_i-1}{z-t_i} + 
\dfrac{\kappa_j+\kappa_k-2}{z-t_j}, \label{eqn:w1} \\
w_2(z) &=& \dfrac{p}{z-q} - \dfrac{L}{z-t_i} 
+ \dfrac{M}{z-t_j} + \dfrac{N}{(z-t_j)^2}, \label{eqn:w2}   
\end{eqnarray}   
where $L$, $M$ and $N$ are given by  
\[
\begin{array}{rcl} 
t_{ij}^2 L &=& q_i q_j^2 p^2 
-\{(\kappa_i-1)q_j + (\kappa_j+\kappa_k) q_i \} q_j p 
+ \kappa_0(\kappa_0+\kappa_4) q_i,  \\
\vspace{-0.3cm} & & \\
t_{ij}^2 M &=&  q_i q_j^2 p^2 
-\{q_i^2+(\kappa_j+\kappa_k-2)q_iq_j+\kappa_i q_j^2\} p 
+ \kappa_0(\kappa_0+\kappa_4) q_i,  \\
\vspace{-0.3cm} & & \\
t_{ij} N &=&  q_i q_j^2 p^2 
-\{(\kappa_j+\kappa_k-1)q_i + \kappa_i q_j\} q_j p 
+ \kappa_0(\kappa_0+\kappa_4) q_j. 
\end{array}  
\]
\end{lem}
{\it Proof}. 
Formula (\ref{eqn:w1}) readily follows from  
the first formula of (\ref{eqn:coeff3}) by letting 
$t_k \to t_j$.  
Next we shall show formula (\ref{eqn:w2}). 
In view of (\ref{eqn:hamiltonian3}), we notice that 
there is a relation $h_i+h_j+h_k = p$ and that $L$, $M$ 
and $N$ are defined in such a manner that 
\[
L = \lim_{t_k \to t_j} h_i, \qquad 
M = L - p, \qquad 
N = - \lim_{t_k \to t_j} t_{jk} h_j.  
\]
In particular, we have $h_j+h_k = p - h_i \to p - L = -M$ 
as $t_k \to t_j$.    
Hence we have  
\[ 
\begin{array}{rcl}
\dfrac{h_i}{z-t_i} &\to& \phantom{-}\dfrac{L}{z-t_i},  \\
\vspace{-0.2cm} & &  \\
\dfrac{h_j}{z-t_j} + \dfrac{h_k}{z-t_k} = 
\dfrac{h_j+h_k}{z-t_k} + \dfrac{t_{jk} h_j}{(z-t_j)(z-t_k)} 
&\to& -\dfrac{M}{z-t_j}-\dfrac{N}{(z-t_j)^2},    
\end{array}
\]   
as $t_k \to t_j$.    
Therefore the second formula of (\ref{eqn:coeff3}) 
leads to (\ref{eqn:w2}). 
The proof is complete. \hfill $\Box$ 
\par\medskip 
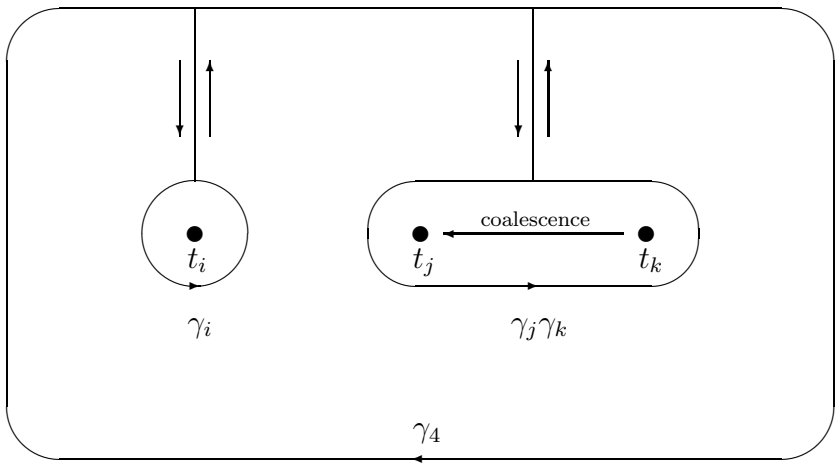
\begin{figure}[t] 
\begin{center}
\setlength{\unitlength}{1.0mm}
\begin{picture}(110,60)(-5,0)
\put(50,30){\oval(110,60)} 
\put(20,30){\circle*{2}}
\put(50,30){\circle*{2}} 
\put(80,30){\circle*{2}}
\put(20,30){\circle{14}}
\put(65,30){\oval(44,14)} 
\put(20.5,23){\vector(1,0){0.2}} 
\put(65.5,23){\vector(1,0){0.2}}
\put(49,0){\vector(-1,0){0.2}}
\put(20,37){\line(0,1){23}}
\put(65,37){\line(0,1){23}}
\put(18,53){\vector(0,-1){10}}
\put(22,43){\vector(0,1){10}}
\put(63,53){\vector(0,-1){10}}
\put(67,43){\vector(0,1){10}}
\put(77,30){\vector(-1,0){24}} 
\put(19,25.5){$t_i$} 
\put(49,25.5){$t_j$}
\put(79,25.5){$t_k$}   
\put(19,17){$\gamma_i$} 
\put(62,17){$\gamma_j\gamma_k$}
\put(49,3){$\gamma_4$}  
\put(58,31){$\ss\mathrm{coalescence}$} 
\end{picture} 
\end{center}
\caption{Coalescence of singular points $t_k \to t_j$} 
\label{fig:coals} 
\end{figure} 
Note that the local exponents at $z = t_i$, $q$, $\infty$ 
and the apparentness of $q$ are preserved in the process 
of coalescence (\ref{eqn:fuchs3}) $\to$ (\ref{eqn:coalfuchs}).     
Let $\lambda_1$ and $\lambda_2$ be the exponents of 
(\ref{eqn:coalfuchs}) at the coalescent singularity $z = t_j$, 
which are the roots of the quadratic equation 
$\lambda^2 - (\kappa_j+\kappa_k-1) \lambda + N = 0$. 
We are interested in its discriminant, that is, the 
squared difference between $\lambda_1$ and $\lambda_2$:   
\begin{equation} \label{eqn:discrim}
\D = (\lambda_1-\lambda_2)^2 = (\kappa_j+\kappa_k-1)^2 - 4 N. 
\end{equation} 
Then the following lemma will play a key role in establishing 
our main theorem.  
\begin{lem}[Key Lemma] \label{lem:key} 
The discriminant $\D$ is $G$-invariant. 
\end{lem}
{\it Proof}. 
Since $t_{ij}$ is $G$-invariant, it is sufficient to show 
the $G$-invariance of $D = -t_{ij} \D$.   
Using $t_{ij} = q_j-q_i$ we have 
$D = (\kappa_j+\kappa_k-1)^2(q_i-q_j) + 4 t_{ij}N$. 
Since $t_{ij}N$ is a polynomial of $(q_i,q_j,p,\kappa)$, 
so is $D$; explicitly, 
$D = D(q_i,q_j,p,\kappa)$ is given by   
\begin{equation} \label{eqn:D}
\begin{array}{rcl}
D(q_i,q_j,p,\kappa) &=& (\kappa_j+\kappa_k-1)^2 (q_i-q_j) \\
\vspace{-0.3cm} & & \\
 &+& 4 q_j [q_iq_jp^2-\{(\kappa_j + \kappa_k-1)q_i + 
\kappa_i q_j\}p+\kappa_0(\kappa_0+\kappa_4)]. 
\end{array} 
\end{equation}
For the moment, we think of $(q_i,q_j,p,\kappa)$ as  
independent variables, namely, we do not assume 
the relations (\ref{eqn:invariant}) and 
(\ref{eqn:fuchsrel}), and put $q_k = q_j$. 
Then a direct check shows that  
\begin{eqnarray*}
s_0(D)-D &=& -4 \kappa_0 (2 q_jp + \kappa_0) \eta, \\ 
s_i(D)-D &=& \- 4 \kappa_i q_j \eta, \\
s_j(D)-D &=& \- 4 \kappa_j q_j \eta, \\
s_k(D)-D &=& \- 4 \kappa_k q_j \eta, \\
s_4(D)-D &=& \- 0,  
\end{eqnarray*} 
where $\eta = 2\kappa_0+\kappa_1+\kappa_2+\kappa_3+\kappa_4-1$. 
Here we must take $q_k = q_j$ and Remark \ref{rem:s4} into 
account in the evaluations of $s_k(D) - D$ and $s_4(D) - D$, 
respectively.     
Hence, under Fuchs' relation (\ref{eqn:fuchsrel}), 
$D$ is an invariant of $s_0$, $s_1$, $s_2$, $s_3$, $s_4$ and so 
is an invariant of $G$. 
The proof is complete.  \hfill $\Box$ 
\par\medskip 
The normalization procedure (\ref{eqn:fuchs3}) $\to$ 
(\ref{eqn:Fuchs3}) passes through the coalescence process, 
leading to a parallel normalization of equation 
(\ref{eqn:coalfuchs}). 
The corresponding gauge transformation is 
\begin{equation} \label{eqn:gauge2} 
f = \psi F \qquad \mathrm{with} \qquad 
\psi = (z-q) (z-t_i)^{\kappa_i/2}(z-t_j)^{(\kappa_j+\kappa_k)/2},        
\end{equation} 
by which equation (\ref{eqn:coalfuchs}) is normalized into a 
Fuchsian equation 
\begin{equation} \label{eqn:coalfuchs2}
\dfrac{d^2F}{dz^2} - W_1(z) \dfrac{dF}{dz} + W_2(z)F = 0,       
\end{equation}   
having singular points at $z = t_i$, $t_j$, $q$, $\infty$, 
such that  
\begin{enumerate} 
\item[(D1)] $z = t_i$ is a regular singular point with exponents 
$\pm\kappa_i/2$;  
\vspace{-0.2cm}
\item[(D2)] $z = q$ is an apparent singular 
point with exponents $\pm1$;  
\vspace{-0.2cm}
\item[(D3)] $z = \infty$ is a regular singular point with 
exponents $(3\pm\kappa_4)/2$;   
\vspace{-0.2cm}
\item[(D4)] $z = t_j$ is a regular singular point with exponents 
$(-1 \pm \sqrt{\D})/2$.    
\end{enumerate} 
\begin{rem} \label{rem:d4}
It follows from (D4) that the trace of the monodromy matrix 
around $z = t_j$ is $- 2 \cos \pi \sqrt{\D}$. 
\end{rem}
\par 
We proceed to taking the isomonodromic deformation into account 
and shall complete the proof of Theorem \ref{thm:main}, leaving 
a certain technical issue to the next section. 
\par\medskip\noindent  
{\it Proof of Theorem $\ref{thm:main}$}.  
The notation of \S\ref{sec:main} will be retained in the 
subsequent discussions, except that the transformation $r$ in 
diagram (\ref{eqn:CD}) will be written $s$ upon identifying $r$ 
and $s$.  
Establishing Theorem \ref{thm:main} amounts to showing that 
\begin{equation} \label{eqn:main} 
s(x_i) = x_i \qquad (i = 1,2,3),      
\end{equation} 
where $x_i = x_i(q,p,t,\kappa)$, which was defined 
in (\ref{eqn:xi}), is a holomorphic function of 
$(q,p,t,\kappa)$ and $s(x_i)$ is understood to be 
$x_i(s(q),s(p),s(t),s(\kappa))$.  
It is sufficient to show (\ref{eqn:main}) for the generators 
$s = s_j$ of the B\"acklund transformation group $G$. 
But this claim for $j = 1, 2, 3, 4$ is trivial, since $s_j$ 
in this case is a gauge transformation as seen in 
\S\ref{sec:gauge}, which does not change the monodromy 
representation and hence the value of $x_i$.  
So the substantial part of the proof is to verify the claim 
for $j = 0$, though the reasoning below will carry over for 
$j = 1,2,3,4$ as well.   
Moreover, we have only to establish it for generic values of 
$\kappa$; the validity of (\ref{eqn:main}) for each $\kappa$ 
in  a dense open subset leads to the validity for every $\kappa$ 
due to the continuous dependence of 
$x_i = x_i(t,q,p,\kappa)$ upon $\kappa$.  
As such a dense open condition on $\kappa$, we employ 
\begin{equation} \label{eqn:generic} 
1-\kappa_j-\kappa_k \in \C\setminus\mathbb{R} 
\qquad (1 \le j < k \le 3),   
\end{equation} 
for a certain technical reason to be explained in 
Lemma \ref{lem:accumulation}.  
\par
The idea of {\it coalescence along isomonodromic deformation} 
proceeds as follows.  
Given a point $(q,p) \in (\C \setminus \{t_1,t_2,t_3\}) 
\times \C$, let $(q(t), p(t))$ be the solution trajectory to 
the Hamiltonian system (\ref{eqn:HS3}) starting from the 
initial point $(q,p)$. 
Since (\ref{eqn:HS3}) describes the isomonodromic deformation 
of Fuchsian equations (\ref{eqn:Fuchs3}), their monodromy 
matrix along the loop $\gamma_j\gamma_k$ 
(see Figure \ref{fig:coals}),  
\[
M_jM_k = (M_jM_k)(q(t),p(t),t,\kappa),   
\]
is independent of $t$ on the trajectory $(q(t),p(t))$.   
Assume that the trajectory admits an accumulation point 
$(\bar{q},\bar{p})$ as $t_k$ tends to $t_j$ along a curve. 
As $t_k$ tends to $t_j$ in such a manner, the Fuchsian 
equation (\ref{eqn:Fuchs3}) degenerates into 
(\ref{eqn:coalfuchs2}), where $(q,p)$ in (\ref{eqn:coalfuchs2})  
should be replaced by $(\bar{q},\bar{p})$.    
In this process, the role of the matrix $M_jM_k$ 
changes from being the global monodromy matrix of 
(\ref{eqn:Fuchs3}) along the loop $\gamma_j\gamma_k$ into 
being the local monodromy matrix of (\ref{eqn:coalfuchs2}) at 
the coalescent singular point $z = t_j$. 
The former role gives $x_i = \Tr(M_jM_k)$, while the latter 
yields $\Tr(M_jM_k) = - 2 \cos \pi \sqrt{\D}$ by virtue of 
Remark \ref{rem:d4}.  
Therefore we have  
\begin{equation} \label{eqn:xiD}
x_i = - 2 \cos \pi \sqrt{\D}. 
\end{equation}     
We can apply the same procedure with the initial point $(q,p)$ 
replaced by $(s(q),s(p))$.  
Since $(s(q(t)), s(p(t)))$ has an accumulation point 
$(s(\bar{q}), s(\bar{p}))$, the same reasoning as above yields   
\begin{equation} \label{eqn:sxiD}
s(x_i) = - 2 \cos \pi \sqrt{s(\D)}. 
\end{equation}     
Thanks to Lemma \ref{lem:key}, which asserts that $\D$ is 
an $s$-invariant, (\ref{eqn:xiD}) and (\ref{eqn:sxiD}) lead 
to the desired equality (\ref{eqn:main}). 
This completes the proof of Theorem \ref{thm:main}.     
\par 
For the true end of the proof, however, there is an extra 
issue yet to be argued, namely, the existence of the 
accumulation point $(\bar{q}, \bar{p})$. 
More precisely, the existence must be assured in such a 
manner that $(\bar{q},\bar{p})$ is located in a general 
position.  
To understand what this means, we should notice that the 
arguments leading to 
Lemmas \ref{lem:coales} and \ref{lem:key} are valid only 
under the condition that the apparent singular point $q$ is  
different from the other singular points $t_i$, $t_j$, 
$\infty$ and that $p$ is a finite complex number.    
In the current situation, this condition should be applied 
to $(q,p) = (\bar{q}, \bar{p})$ as well as to $(q,p) = 
(s(\bar{q}), s(\bar{p}))$.  
For the generators $s = s_0$, $s_1$, $s_2$, $s_3$, $s_4$ 
of $G$, formula (\ref{eqn:si}) implies that the condition 
applied to them is expressed as       
\begin{equation} \label{eqn:accumulation}
\bar{q} \in \C \setminus \{t_i,t_j\}, \qquad 
\bar{p} \in \C \setminus \{0\}, \qquad 
\bar{q} + \dfrac{\kappa_0}{\bar{p}}\in \C \setminus 
\{t_i,t_j\}.    
\end{equation}  
Thus an accumulation point $(\bar{q}, \bar{p})$ is said 
to be in a general position if it satisfies condition  
(\ref{eqn:accumulation}). 
Intuitively the existence of such a point is quite likely,  
but logically nontrivial. 
It is at this stage that the generic condition 
(\ref{eqn:generic}) on $\kappa$ is used to provide 
the following lemma.                      
\begin{lem} \label{lem:accumulation} 
Suppose that 
$\kappa = (\kappa_0,\kappa_1,\kappa_2,\kappa_3,\kappa_4) 
\in \K$ satisfies condition $(\ref{eqn:generic})$. 
Then there exists an open subset $V$ of 
$(\C\setminus\{t_1,t_2,t_3\}) \times \C$ such that 
for every $(q,p) \in V$, the solution trajectory $(q(t), p(t))$ 
to the Hamiltonian system $(\ref{eqn:HS3})$ starting from the 
initial point $(q, p)$ admits an accumulation point 
$(\bar{q},\bar{p})$ satisfying condition 
$(\ref{eqn:accumulation})$ as $t_k$ tends to $t_j$ along 
a curve.      
\end{lem}
\par 
Lemma \ref{lem:accumulation} guarantees the existence 
of an accumulation point $(\bar{q},\bar{p})$ in a general 
position only when the initial point $(q,p)$ belongs to an 
open subset $V$. 
But this is enough for proving (\ref{eqn:main}). 
Indeed, by Lemma \ref{lem:accumulation}, equation 
(\ref{eqn:main}) is valid at least for $(q,p) \in V$. 
Then the unicity theorem for holomorphic functions 
implies that it remains valid for every $(q,p)$, 
since $x_i = x_i(q,p,t,\kappa)$ is holomorphic in 
$(q,p)$ and the space $\E_t(\kappa) \cong 
(\C \setminus \{t_1,t_2,t_3\}) \times \C$ is connected. 
To establish Theorem \ref{thm:main}, it only remains to 
prove Lemma \ref{lem:accumulation}.  
This final task will be done in \S\ref{sec:accumulation}. 
\hfill $\Box$     
\section{Accumulation Points in a General Position} 
\label{sec:accumulation} 
We shall establish Lemma \ref{lem:accumulation} for $j = 1$ 
and $k = 3$, namely, for the coalescence process $t_3 \to t_1$; 
due to the symmetry in $t_i$, $t_j$, $t_k$, the other cases 
can be treated in a similar manner.  
For this purpose, we may work with the Hamiltonian system 
(\ref{eqn:HS}) with single time variable $x$ instead of 
the system (\ref{eqn:HS3}) with three time variables 
$t = (t_1,t_2,t_3)$.  
In view of (\ref{eqn:cross}), letting $t_3 \to t_1$ means 
$x \to 0$. 
With these remarks, Lemma \ref{lem:accumulation} is reduced to 
the following:         
\begin{lem} \label{lem:accumulation2}  
Suppose that 
$\kappa = (\kappa_0, \kappa_1,\kappa_2,\kappa_3,\kappa_4) \in \K$ 
satisfies a generic condition 
\begin{equation} \label{eqn:eta} 
1-\kappa_1-\kappa_3 \in \C\setminus\mathbb{R}. 
\end{equation}
Then there exists a $2$-parameter family of solutions to the 
Hamiltonian system $(\ref{eqn:HS})$, 
\begin{equation} \label{eqn:twoparam}
(q(x,c), p(x,c)), \qquad c = (c_1, c_2) \in U, 
\end{equation} 
where $U$ is an open subset of $\C^2$, such that 
the following conditions are satisfied$:$ 
\begin{enumerate}
\item The correspondence $c = (c_1,c_2) \mapsto 
(q(x_0,c), p(x_0,c))$  defines a biholomorphic map of $U$ onto 
an open subset $V$ of $(\C\setminus\{0,1\}) \times \C$, 
where $x_0$ is a point in $\C \setminus \{0,1\}$.  
\item The solution $(\ref{eqn:twoparam})$ admits an accumulation 
point $(q(c), p(c))$ as $x \to 0$ such that 
\begin{equation} \label{eqn:accumulation2} 
q(c) \in \C \setminus \{0,1\}, \qquad 
p(c) \in \C \setminus \{0\}, \qquad 
q(c) + \dfrac{\kappa_0}{p(c)} \in \C \setminus \{0,1\}. 
\end{equation} 
\end{enumerate} 
\end{lem} 
\par 
Note that condition (\ref{eqn:accumulation2}) corresponds to 
(\ref{eqn:accumulation}).  
For a proof of Lemma \ref{lem:accumulation2}, we utilize a result 
by Takano \cite{Takano} and Kimura \cite{Kimura}, who established 
a reduction theorem and constructed a $2$-parameter family of 
solutions to $\PVI$ around its fixed singular points. 
We recall their construction in a manner suitable for our purpose. 
Put 
\[
\begin{array}{rcl}
E(r,\rho) &=&  \{\, (x,Q,P) \in \C^3\,:\, |x| < r, \, 
|Q| < \rho, \, |xP| < \rho, \, |QP| < \rho \,\}, \\
\vspace{-0.3cm} & & \\
E(\rho) &=&  \{\, (Q,P) \in \C^2\,:\, |Q| < \rho, \, 
|QP| < \rho \,\}. 
\end{array} 
\] 
Then the following lemma is an easy consequence of 
Theorems 1 and 2 in Takano \cite{Takano}. 
\begin{lem} \label{lem:takano} 
Suppose that $\kappa = 
(\kappa_0,\kappa_1,\kappa_2,\kappa_3,\kappa_4) \in \K$ 
satisfies condition $(\ref{eqn:eta})$. 
Then there exist positive constants $r$, $\rho > 0$ and a unique 
canonical transformation 
\[
q = b(x,Q,P), \qquad p = a(x,Q,P) 
\]
that reduces system $(\ref{eqn:HS})$ into a Hamiltonian system 
with Hamiltonian 
\begin{equation} \label{eqn:HSquad} 
h_0(x,Q,P) = \{(QP)^2 + (1-\kappa_1-\kappa_3)(QP) \}/x,  
\end{equation}
where $b(x,Q,P)$ and $a(x,Q,P)$ are holomorphic functions in 
$E(r,\rho)$ such that 
\begin{equation} \label{eqn:ba}
|b(0,Q,P)-Q| \le M|Q|^2, \qquad |a(0,Q,P) - P| \le M 
\qquad\mathrm{for} \quad (Q,P) \in E(\rho),  
\end{equation} 
with some positive constant M$;$ we may and shall assume 
that $M > 2$. 
\end{lem} 
\par 
The Hamiltonian system with Hamiltonian 
(\ref{eqn:HSquad}) has a first integral $QP$ and is settled as 
\begin{equation} \label{eqn:lambda} 
Q(x,c) = c_1 x^{\lambda}, \quad 
P(x,c) = c_2 x^{-\lambda}  \qquad \mathrm{with} \quad 
\lambda = \lambda(c) = 1-\kappa_1-\kappa_3 + 2 c_1c_2, 
\end{equation} 
where $c = (c_1,c_2) \in \C^2$ is an arbitrary constant. 
Then Lemma \ref{lem:takano} asserts that 
\begin{equation} \label{eqn:qp}
q(x,c) = b(x, Q(x,c), P(x,c)), \qquad 
p(x,c) = a(x, Q(x,c), P(x,c)), 
\end{equation} 
yield a solution to system (\ref{eqn:HS}), provided that 
$(x,Q(x,c),P(x,c)) \in E(r,\rho)$, namely,  
\[
0 < |x| < r, \qquad |Q(x,c)| < \rho, \qquad |x P(x,c)| < \rho, 
\qquad |c_1c_2| < \rho. 
\]
In terms of $\arg x$ and $\log |x|$, the second and third 
conditions are expressed as in Table \ref{tab:domain},  
the exhibition of which is divided into five cases according 
to the values of $\Re\,\lambda$. 
\par    
\begin{table}[t] 
\begin{center}
\begin{tabular}{|c||c|c|c|}
\hline
& & & \\
case & $\lambda$ & $\arg x$ & $\log|x|$ \\
& & & \\
\hline
\hline
& & & \\
1 & $\Re\, \lambda > 1$ & 
$(\Im\,\lambda)(\arg x) < L_0$ &  
$L_2(\arg x) < \log |x| < L_1(\arg x)$ \\
& & & \\
\hline
& & & \\
2 & $\Re\,\lambda = 1$ & 
$(\Im\,\lambda)(\arg x) < L_0$ & 
$\log |x| < L_1(\arg x)$ \\
& & & \\
\hline
& & & \\
3 & $0 < \Re\,\lambda < 1$ & no constraints & 
$\log |x| < \min \{L_1(\arg x), \, L_2(\arg x) \}$ \\
& & & \\
\hline
& & & \\
4 & $\Re\,\lambda = 0$ & 
$(\Im\,\lambda)(\arg x) > L_0$ & 
$\log|x| < L_2(\arg x)$ \\
& & & \\
\hline
& & & \\
5 & $\Re\, \lambda < 0$ & 
$(\Im\,\lambda)(\arg x) > L_0$ & 
$L_1(\arg x) < \log |x| < L_2(\arg x)$ \\
& & & \\ 
\hline
\end{tabular}
\end{center}
\[
\begin{array}{rcl}
L_0 &=& (\Re\,\lambda)\log(\rho/|c_2|) + 
(\Re\,\lambda-1)\log(\rho/|c_1|), \\
\vspace{-0.2cm} & & \\
L_1(\alpha) &=& \dfrac{(\Im\,\lambda)\alpha + 
\log(\rho/|c_1|)}{\Re\,\lambda}, \\
\vspace{-0.2cm} & & \\
L_2(\alpha) &=& \dfrac{(\Im\,\lambda)\alpha - 
\log(\rho/|c_2|)}{\Re\,\lambda-1}.
\end{array} 
\] 
\caption{The conditions $|Q(x,c)| < \rho$ and $|x P(x,c)| < \rho$} 
\label{tab:domain} 
\end{table} 
\medskip\noindent 
{\it Proof of Lemma $\ref{lem:accumulation2}$}.  
By assumption (\ref{eqn:eta}), we can choose a number $\rho_0$ 
so that  
\begin{equation} \label{eqn:rho0} 
0 < \rho_0 < \min
\{\rho, \, |\Im(1-\kappa_1-\kappa_3)|/2, \, 1 \}.   
\end{equation} 
Out of later necessity, if $\kappa_0 \neq 0$, we take $\rho_0$ 
so as to satisfy an additional condition  
\begin{equation} \label{eqn:rho00}
\rho_0 < \dfrac{|\kappa_0|}{2},    
\end{equation}
which should be neglected if $\kappa_0 = 0$.  
As the open subset $U$ mentioned in (\ref{eqn:twoparam}), 
we put   
\begin{equation} \label{eqn:U}
U = \{\, c = (c_1,c_2) \in \C^2\,:\, 0 < |c_1c_2| < \rho_0 \,\}. 
\end{equation}  
For each $c \in U$, let $D(c)$ be the set of those points on the 
universal covering of $0 < |x| < r$ which satisfy the conditions 
in Table \ref{tab:domain}.  
By (\ref{eqn:lambda}), (\ref{eqn:rho0}) and (\ref{eqn:U}), one 
has $\Im\,\lambda(c) \neq 0$, which implies that $D(c)$ is a 
nonempty domain that contains a curve tending to the 
origin $x = 0$. 
Then the solution $(q(x,c), p(x,c))$ in (\ref{eqn:qp}) makes 
sense on the domain $D(c)$. 
\par
\begin{figure}[t] 
\[
\begin{array}{ccc}
\quad \includegraphics*[width=4cm,height=4cm]{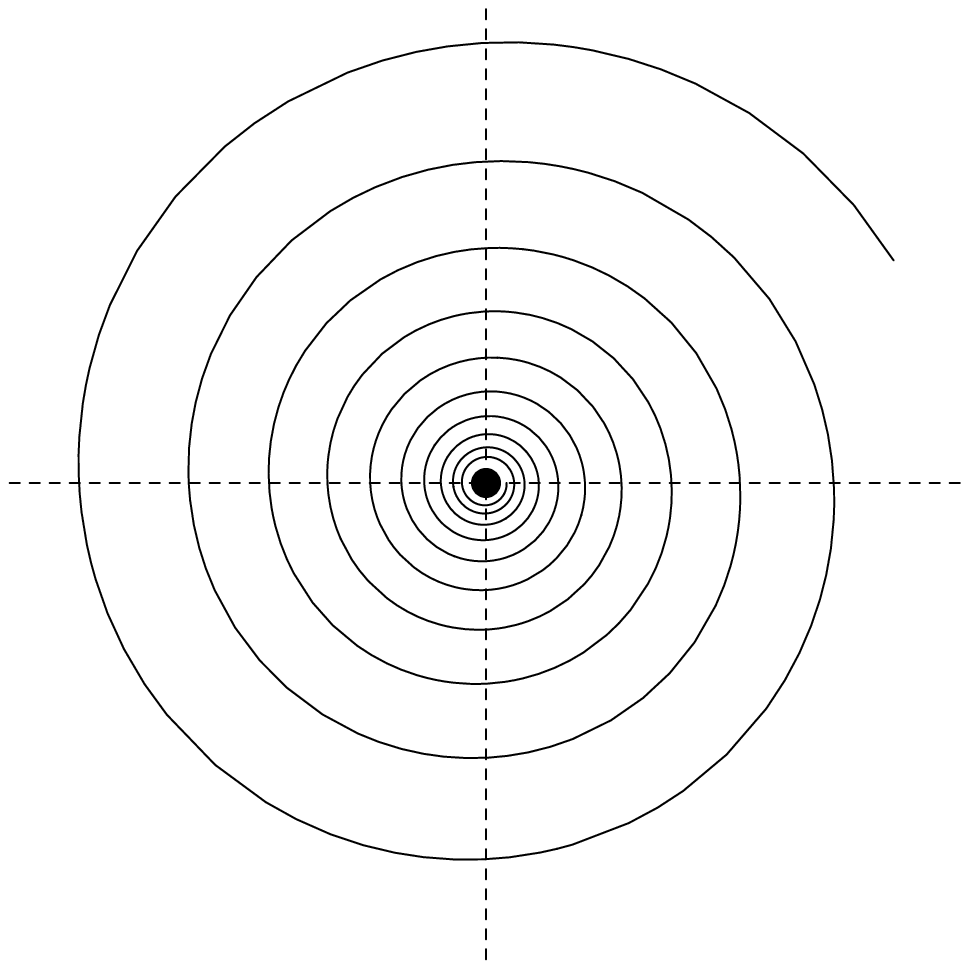} \quad & 
\quad \includegraphics*[width=4cm,height=4cm]{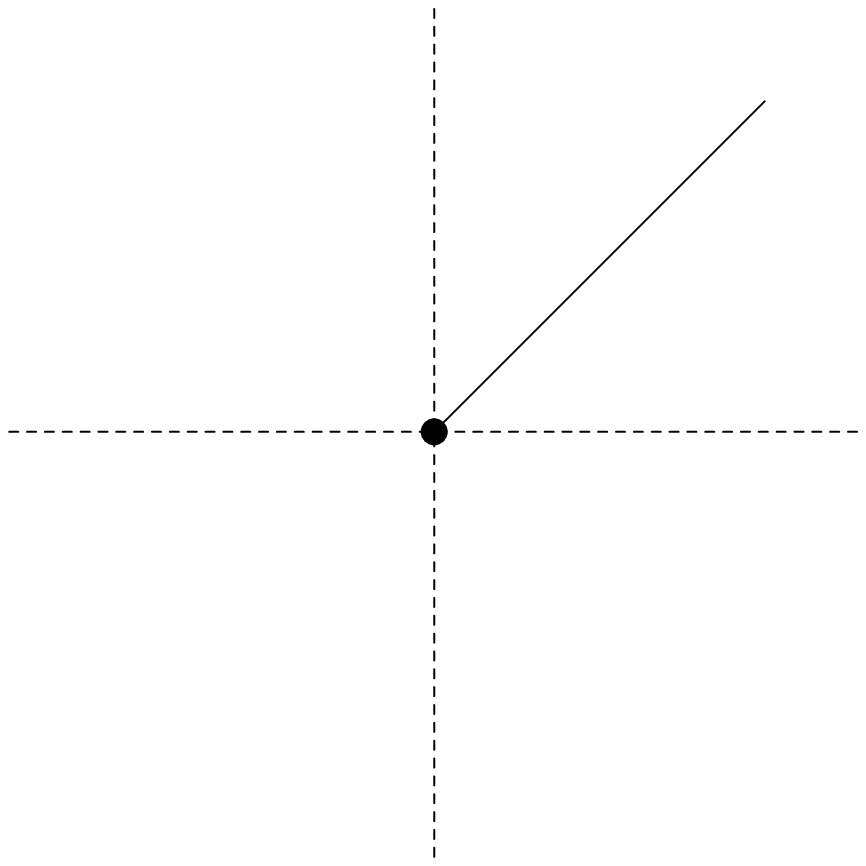} \quad &
\quad \includegraphics*[width=4cm,height=4cm]{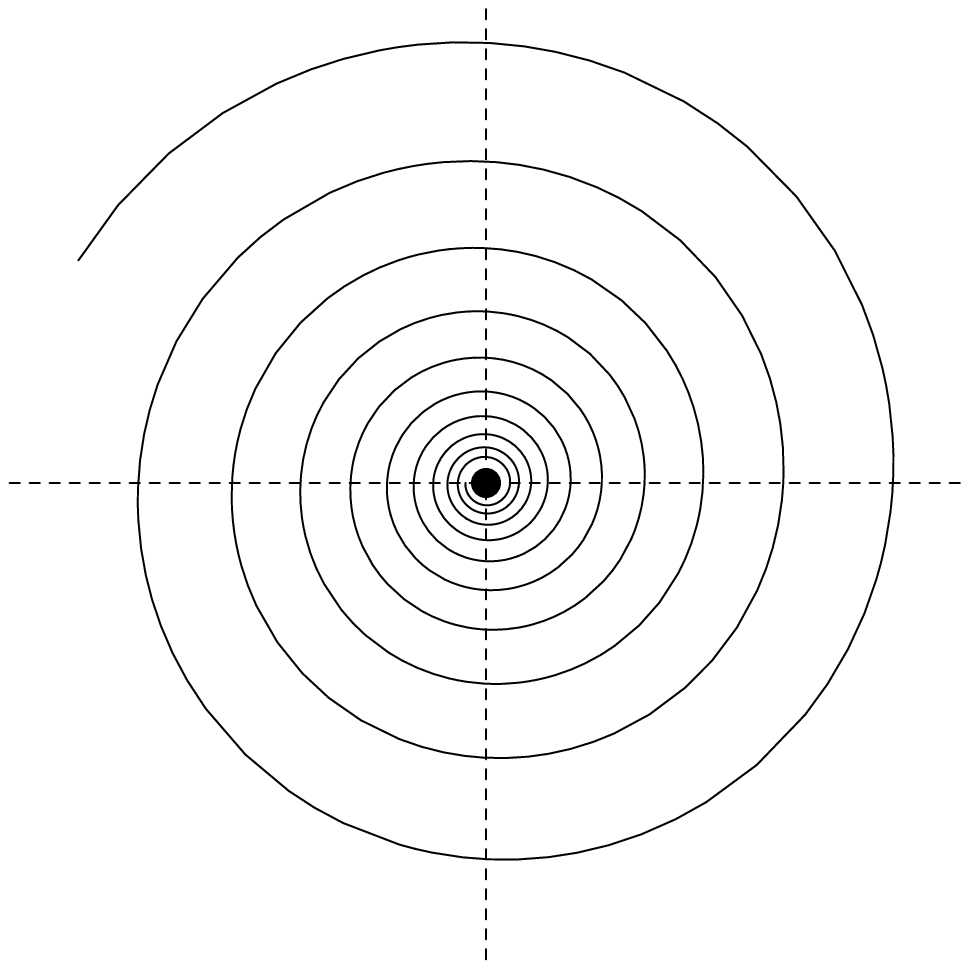} \quad \\
      &    &     \\
\quad (\Re\,\lambda)(\Im\,\lambda ) < 0 & 
\Re\,\lambda = 0 & 
\quad (\Re\,\lambda)(\Im\,\lambda ) > 0 
\end{array} 
\]
\caption{The curve $\gamma$: a spiral or a line segment} 
\label{fig:spiral} 
\end{figure}

 
For each $c \in U$, we shall show that $(q(x,c), p(x,c))$ 
admits an accumulation point $(q(c),p(c))$ satisfying 
condition (\ref{eqn:accumulation2}) as $x \to 0$ along a 
curve in $D(c)$.   
Let $\mu$ be a number such that  
\begin{equation} \label{eqn:mu} 
0 < \mu < \dfrac{|c_1c_2|}{M(|\kappa_0|+8)}.     
\end{equation}
and consider a curve in the $x$-plane defined by  
\[
\gamma = \left\{\,x\,:\, (\Re\,\lambda) \log|x| - 
(\Im\,\lambda)(\arg x) = \log\dfrac{\mu}{|c_1|}, \,\, 
|x|< \dfrac{\mu\rho}{|c_1c_2|} \,\right\}    
\]
Note that (\ref{eqn:rho0}), (\ref{eqn:U}), (\ref{eqn:mu}) 
and $M > 2$ imply $\mu < \rho$. 
Along the curve $\gamma$ one has  
\begin{equation} \label{eqn:absolute}
|Q(x,c)|= \mu, \qquad |P(x,c)| = \dfrac{|c_1c_2|}{\mu}.
\end{equation} 
In particular, $|Q(x,c)| < \rho$ and $|xP(x,c)| < \rho$, 
and hence the curve $\gamma$ lies in the domain $D(c)$. 
The shape of $\gamma$ is indicated in Figure \ref{fig:spiral}; 
if $\Re\,\lambda \neq 0$, it is a spiral curve winding around 
and tending to the origin, while if $\Re\,\lambda = 0$, it 
is a line segment terminating at the origin. 
\par  
In view of (\ref{eqn:absolute}), as $x$ tends to the origin 
along the curve $\gamma$, there exists an accumulation 
point $(Q(c),P(c)) \in E(\rho)$ of $(Q(x,c),P(x,c))$ such 
that $Q(c)P(c) = c_1c_2$ and 
\begin{equation} \label{eqn:absolute2}
|Q(c)|= \mu, \qquad |P(c)| = \dfrac{|c_1c_2|}{\mu}. 
\end{equation} 
Accordingly, as $x \to 0$ along $\gamma$, $(q(x,c),p(x,c))$ 
has an accumulation point $(q(c),p(c))$ with  
\[
q(c) = b(0,Q(c),P(c)), \qquad p(c) = a(0,Q(c),P(c)).  
\] 
We shall show that $(q(c),p(c))$ satisfies the desired 
property (\ref{eqn:accumulation2}).  
Hereafter $(Q(c), P(c))$ and $(q(c), p(c))$ will be  
abbreviated to $(Q,P)$ and $(q,p)$ respectively. 
The arguments below are rather technical, but the idea 
itself is very simple: If $\mu > 0$ is sufficiently small, 
then (\ref{eqn:absolute2}) implies that $Q$ is sufficiently 
small but not zero and $P = c_1c_2/Q$ is sufficiently 
large but finite.  
Now (\ref{eqn:ba}) means that $(q,p)$ equals $(Q,P)$ 
up to a first order error term.    
Then one can show that $q$ is small but not zero, 
$p$ is large but finite and also $q + \kappa_0/p$ is 
small but not zero.  
\par  
We verify the first and second conditions of 
(\ref{eqn:accumulation2}). 
Note that (\ref{eqn:rho0}), (\ref{eqn:U}) and 
(\ref{eqn:mu}) yield  
\[  
M < \dfrac{|c_1c_2|}{2\mu} < \dfrac{1}{2\mu}, \qquad 
\mu < \dfrac{1}{2}. 
\]
Then (\ref{eqn:ba}) and (\ref{eqn:absolute2}), 
combined with these inequalities, lead to   
\[ 
\begin{array}{rcccccccccl}
|q| &\le& |Q|+|q-Q| &\le& |Q|+M|Q|^2 &=& (1+M \mu)\mu 
&<& \dfrac{3\mu}{2} &<& \dfrac{3}{4}, \\
\vspace{-0.2cm} &&&&&&&&&& \\
|q| &\ge& |Q|-|q-Q| &\ge& |Q|-M|Q|^2 &=& (1-M \mu)\mu  
&>& \dfrac{\mu}{2} &>& 0,  \\
\vspace{-0.2cm} &&&&&&&&&& \\ 
|p| &\ge& |P|-|p-P| &\ge& |P|-M &=& 
\dfrac{|c_1c_2|}{\mu} - M &>& 
\dfrac{|c_1c_2|}{2 \mu} &>& 0,      
\end{array}
\] 
which verifies the first and second conditions of 
(\ref{eqn:accumulation2}) as desired.  
\par 
We proceed to verify the third condition of 
(\ref{eqn:accumulation2}).  
We may assume $\kappa_0 \neq 0$; for otherwise the 
third condition is the same as the first one. 
We observe that  
\begin{equation} \label{eqn:estimate}
|p| > \dfrac{|c_1c_2|}{2\mu}, \qquad 
|c_1c_2| < \dfrac{|\kappa_0|}{2}, \qquad 
\mu < \dfrac{|c_1c_2|}{2|\kappa_0|}. 
\end{equation} 
Indeed, the first one is already seen in the last 
paragraph, the second one is obtained from (\ref{eqn:rho00}) 
and (\ref{eqn:U}), and the last one follows from 
(\ref{eqn:mu}) and $M > 2$, respectively.   
Let    
\[
R = \left(q + \dfrac{\kappa_0}{p} \right) - 
\left(1 + \dfrac{\kappa_0}{c_1c_2} \right)Q  
= (q-Q) + \left( \dfrac{\kappa_0}{p} - 
\dfrac{\kappa_0}{P} \right) 
= (q-Q) + \dfrac{\kappa_0(P-p)}{pP}. 
\]   
By applying (\ref{eqn:ba}), (\ref{eqn:absolute2}), 
(\ref{eqn:estimate}) and (\ref{eqn:mu}) 
successively, $R$ is estimated as      
\[ 
\begin{array}{rcccl} 
|R| 
&\le& |q-Q| + \dfrac{|\kappa_0||p-P|}{|p||P|} 
&\le& M|Q|^2 + \dfrac{|\kappa_0|M}{|p||P|} \\
\vspace{-0.2cm} & & & & \\
&<& M \mu^2 + \dfrac{2|\kappa_0| M \mu^2}{|c_1c_2|^2} 
&=& \dfrac{M\mu^2(|c_1c_2|^2 +2|\kappa_0|)}{|c_1c_2|^2} \\
\vspace{-0.2cm} & & & & \\    
&<& \dfrac{M \mu^2 |\kappa_0|(|\kappa_0|+8)}{4|c_1c_2|^2} 
&<& \dfrac{|\kappa_0|\mu}{4|c_1c_2|}. 
\end{array} 
\] 
On the other hand, the second inequality in 
(\ref{eqn:estimate}) leads to  
\[
\begin{array}{rcccccl}
\left|1+\dfrac{\kappa_0}{c_1c_2} \right| &\le&  
\dfrac{|\kappa_0|}{|c_1c_2|} + 1 &<& 
\dfrac{|\kappa_0|}{|c_1c_2|} + \dfrac{|\kappa_0|}{2|c_1c_2|} 
&=& \dfrac{3|\kappa_0|}{2|c_1c_2|}, \\  
\vspace{-0.3cm} &&&&&& \\
\left|1+\dfrac{\kappa_0}{c_1c_2} \right| &\ge&  
\dfrac{|\kappa_0|}{|c_1c_2|} - 1 &>& 
\dfrac{|\kappa_0|}{|c_1c_2|} - \dfrac{|\kappa_0|}{2|c_1c_2|} 
&=& \dfrac{|\kappa_0|}{2|c_1c_2|}. \\
\end{array} 
\]
These preliminary estimates, (\ref{eqn:absolute}) and 
the third inequality in (\ref{eqn:estimate}) yield 
\[
\begin{array}{rcl} 
\left|q + \dfrac{\kappa_0}{p} \right| 
&=& 
\left|\left(1+\dfrac{\kappa_0}{c_1c_2}\right) Q + R \right| 
\le \left| 1+\dfrac{\kappa_0}{c_1c_2} \right| |Q| + |R| 
< \dfrac{3|\kappa_0|\mu}{2|c_1c_2|} + 
  \dfrac{|\kappa_0|\mu}{4|c_1c_2|} 
= \dfrac{7|\kappa_0|\mu}{4|c_1c_2|} 
< \dfrac{7}{8}, \\ 
\vspace{-0.2cm} & & \\ 
\left|q + \dfrac{\kappa_0}{p} \right| 
&=&  
\left|\left(1+\dfrac{\kappa_0}{c_1c_2}\right) Q + R \right| 
\ge \left| 1+\dfrac{\kappa_0}{c_1c_2} \right| |Q| - |R| 
> \dfrac{|\kappa_0|\mu}{2|c_1c_2|} -  
  \dfrac{|\kappa_0|\mu}{4|c_1c_2|} 
= \dfrac{|\kappa_0|\mu}{4|c_1c_2|} 
> 0.   
\end{array}
\]  
Hence the third condition of (\ref{eqn:accumulation2}), i.e., 
assertion (2) of Lemma \ref{lem:accumulation2} is verified. 
The proof of assertion (1) is omitted, since 
it is a standard inverse function argument, in which   
the open set $U$ will be replaced by a smaller one, if 
necessary.    
The proof is complete. \hfill $\Box$  
\section{Finding B\"acklund Transformations} 
\label{sec:heuristics} 
The B\"acklund transformations $s_i$ have been found by  
Okamoto \cite{Okamoto2}, Arinkin and Lysenko \cite{AR2}, 
Noumi and Yamada \cite{NY} and others by various methods. 
In any case, $s_1$, $s_2$, $s_3$, $s_4$ are easy to find, 
since they are elementary gauge transformations as 
constructed in \S\ref{sec:gauge}. 
But things are different with the transformation $s_0$; it 
cannot be a gauge transformation of rank two differential 
equations, since it does change local monodromy data.  
Now we wish to propose another way of finding $s_0$ from 
our point of view.  
\par
Our idea is to revisit Lemma \ref{lem:key} together with  
formula (\ref{eqn:D}), which is a key observation in this 
paper. 
It asserts that the discriminant ${\mit\Delta}$ in 
(\ref{eqn:discrim}), or equivalently the function $D$ in 
(\ref{eqn:D}), is $G$-invariant. 
Here we write $D = D(q,p,t,\kappa)$, since $D$ is a 
polynomial of $(q,p,t,\kappa)$ with $t_j = t_k$.  
An essence of the reasoning in \S\ref{sec:coals} is the 
implication that if there is an invariance relation 
\begin{equation} \label{eqn:DD} 
D(Q,P,t,\sigma_0(\kappa)) = D(q,p,t,\kappa), 
\end{equation}
then one has $x_i(Q,P,t,\sigma_0(\kappa)) = x_i(q,p,t,\kappa)$ 
for the global monodromy data $x_i = \Tr(M_jM_k)$.  
Now our principle --- a B\"acklund transformation should 
be a pull-back of the identity transformation on the 
moduli of global monodromy data --- suggests that relation 
(\ref{eqn:DD}) would give us a B\"acklund transformation   
$(q,p) \mapsto (Q,P)$. 
So it must be promising to find out $(Q,P)$ as a function  
of $(q,p)$ satisfying (\ref{eqn:DD}). 
This thought brings our attention to the difference 
\begin{equation} \label{eqn:E} 
E = E(Q,P;q,p;t,\kappa) =  
D(Q,P,t,\sigma_0(\kappa)) - D(q,p,t,\kappa).  
\end{equation}  
Consider $E$ as a polynomial of $(t_i,t_j)$ and let $E_{mn}$ 
denote the coefficient of the term $t_i^mt_j^n$ in $E$. 
If $E \equiv 0$, each coefficient $E_{mn}$ must vanish. 
We especially take a look at $E_{12}$ and $E_{11}$:          
\[
\left\{ 
\begin{array}{rcl}
E_{12} &=& 4(p - P)(p + P),  \\
\vspace{-0.3cm} & & \\
E_{11} &=& 4p(\kappa_j + \kappa_k - 1 - 2qp) 
+ 4P(\kappa_i + \kappa_4 + 2QP). 
\end{array}
\right.  
\] 
As is easily seen, 
the system of equations $E_{12} = E_{11} = 0$ has two 
solutions 
\begin{eqnarray} 
Q &=& q + \dfrac{\kappa_0}{p}, 
\phantom{+\kappa_i+\kappa_4} \qquad\,\,  
P \, = \, \phantom{-}p,  
\label{eqn:sol1} \\
Q &=& q + \dfrac{\kappa_0+\kappa_i+\kappa_4}{p},   
\qquad P \, = \, -p. \label{eqn:sol2} 
\end{eqnarray}  
The first solution (\ref{eqn:sol1}) is none other 
than the B\"acklund transformation $s_0$ we seek.   
In this case we can check that substituting 
(\ref{eqn:sol1}) into (\ref{eqn:E}) yields $E \equiv 0$. 
\par 
On the other hand, the second solution (\ref{eqn:sol2}) 
does not imply $E \equiv 0$ in general. 
Indeed, (\ref{eqn:sol2}) leads to    
$E_{02} = 4(2\kappa_i-\kappa_j-\kappa_k+1)p$, and 
hence $E$ does not vanish unless $\kappa$ satisfies    
\begin{equation} \label{eqn:E02} 
2\kappa_i-\kappa_j-\kappa_k+1 = 0. 
\end{equation} 
Assume that (\ref{eqn:E02}) is the case. 
Then, substituting (\ref{eqn:sol2}) and (\ref{eqn:E02}) 
into (\ref{eqn:E}), we have    
$E = 2\kappa_i(\kappa_4-\kappa_i)(\kappa_4+\kappa_i)/p$. 
Therefore, (\ref{eqn:sol2}) yields $E \equiv 0$ if and 
only if $\kappa$ also satisfies     
\begin{equation} \label{eqn:EE}
\kappa_i(\kappa_4-\kappa_i)(\kappa_4+\kappa_i) = 0. 
\end{equation} 
We wonder whether the transformation (\ref{eqn:sol2})  
has any meaning under the restriction of parameters, 
(\ref{eqn:E02}) and (\ref{eqn:EE}). 
But this point will not be touched in this paper.  
\par  
It is worth considering whether there is any other 
B\"acklund transformation than those are already known. 
As for this question, the following proposition shows that 
the B\"acklund transformations are exhausted by the known 
ones, even if the class is enlarged to the analytic category.     
\begin{prop}  \label{prop:backlund} 
Let $\pi : \M_t \to \K$ be the family of moduli spaces  
of stable parabolic bundles with connections 
in $\S\ref{sec:RH}$.   
Let $\sigma$ be an analytic automorphism of $\K$ such 
that $\theta(\kappa) = \theta(\sigma(\kappa))$ for   
$\kappa \in \K$ and let $s$ be a bimeromorphic 
automorphism of $\M_t$ such that $\sigma\pi = \pi s$.  
Assume that for general values of $\kappa \in \K$, the 
analytic isomorphisms 
$s_{\kappa} : \M_t(\kappa) \to \M_t(\sigma(\kappa))$ induce 
the identity on $\mathcal{S}(\theta(\kappa))$ via the 
Riemann-Hilbert correspondence $(\ref{eqn:RH3})$. 
Then $\sigma \in W(\mathrm{D}_4^{(1)})$ and $s$ is a 
known B\"acklund transformation, the unique lift 
of $\sigma$.    
\end{prop} 
{\it Proof}. 
>From the invariant-theoretical argument in Terajima 
\cite{Terajima}, it is not difficult to see that the analytic 
quotient $\K/W(\mathrm{D}_4^{(1)})$ is biholomorphic to 
the complex $4$-space $\C^4_{\theta}$ with coordinates 
$\theta = (\theta_1,\theta_2,\theta_3,\theta_4)$. 
In particular, distinct $W(\mathrm{D}_4^{(1)})$-orbits 
in $\K$ have distinct values of $\theta$.   
This observation shows that any analytic automorphism 
$\sigma$ of $\K$ such that 
$\theta(\kappa) = \theta(\sigma(\kappa))$ for 
$\kappa \in \K$ is necessarily an element of 
$W(\mathrm{D}_4^{(1)})$. 
Then clearly the transformation $s$ is obtained as the 
unique lift of $\sigma$ relative to the Riemann-Hilbert 
correspondence.    
\hfill $\Box$       
\par\medskip  
We conclude this paper by putting some questions to 
ourselves.  
We were able to characterize the B\"acklund 
transformations in a natural manner in terms of 
the Riemann-Hilbert correspondence based on the  
following nice observation:   
The difference of the two local exponents at a coalescent 
regular singular point happens to be an invariant of 
the B\"acklund transformations (Lemma \ref{lem:key}).  
Does this phenomenon occur just by chance or more 
universally?  
Do similar phenomena occur for other Painlev\'e equations 
than $\PVI$ or for Garnier systems? 
If so, do they help us find B\"acklund transformations 
for those equations?      
\par\medskip\noindent 
{\bf Acknowledgments}. 
The authors would like to thank Kyoichi Takano for 
helpful discussions.  


\begin{thebibliography}{99}
%
\bibitem{AL1}{Arinkin, D. and Lysenko, S.},
{\sl On the moduli of ${\rm SL}(2)$-bundles with connections on 
$\P^1\backslash \{x_1, \cdots, x_4 \}$ }.
Internat. Math. Res. Notices 1997, no. 19, 983--999.  
%
\bibitem{AR2}{Arinkin, D. and Lysenko, S.}, 
{\sl Isomorphisms between moduli spaces of $SL(2)$-bundles 
with connections on $\P^1 \setminus \{x_1,\dots,x_4\}$}, 
Math. Res. Lett. {\bf 4} (1997), 181--190. 
%
\bibitem{CM}{Conte, R. and Musette, M.}, 
{\sl First-degree birational transformations of the Painlev\'e 
equations and their contiguity relations}, 
Jour. Phys., A: Math. General {\bf 34} (2001), 10507--10522. 
%
\bibitem{FI}{Fokas, A.S. and Its, A.R.}, 
{\sl The isomonodromy method and the Painlev\'e equations}, 
in ``Important development in soliton theory", Springer Ser. 
Nonlinear Dynamics, Springer-Verlag, Berlin, 1993, 99--122. 
%
\bibitem{FY}{Fokas, A.S. and Yortsos, Y.C.}, 
{\sl The transformation properties of the sixth Painlev\'e 
equation and one-parameter families of solutions}, 
Lett. Nuovo Cimento {\bf 30} (1981), 539--544. 
%
\bibitem{Fuchs}{Fuchs, R.}, 
{\sl \"Uber lineare homogene Differentialgleichungen zweiter 
Ordnung mit drei im Endlichen gelegenen wesentlich singul\"aren 
Stellen}, Math. Ann. {\bf 63} (1907), 301--321. 
%
\bibitem{GLS}{Gromak, V.I., Laine, I. and Shimomura, S.}, 
{\sl Painlev\'e differential equations in the complex plane}, 
Studies in Math. {\bf 28}, de Gruyter, Berlin, 2002. 
%
\bibitem{IIS}{Inaba, M., Iwasaki, K. and Saito, M.-H.}, 
{\sl Moduli of stable parabolic connections, Riemann-Hilbert 
correspondence and geometry of Painlev\'e equations of type VI, 
Part I}, Preprint, 2003. 
%
\bibitem{Iwasaki1}{Iwasaki, K.}, 
{\sl Moduli and deformation for Fuchsian projective connections 
on a Riemann surface}, J. Fac. Sci. Univ. Tokyo., Sect. IA, Math. 
{\bf 38} (1991), 431--531.  
%
\bibitem{Iwasaki2}{Iwasaki, K.}, 
{\sl Fuchsian moduli on Riemann surfaces --- its Poisson structure 
and Poincar\'e-Lefschetz duality}, Pacific J. Math. {\bf 155} 
(1992), 319--340. 
%
\bibitem{Iwasaki3}{Iwasaki, K.}, 
{\sl A modular group action on cubic surfaces and the 
monodromy of the Painlev\'e VI equation}, Proc. Japan 
Acad., {\bf 78}, Ser. A (2002), 131--135. 
%
\bibitem{Iwasaki4}{Iwasaki, K.}, 
{\sl An area-preserving action of the modular group on 
cubic surfaces and the Painlev\'e VI equation}, 
to appear in Comm. Math. Phys.   
%
\bibitem{JMU}{Jimbo, M., Miwa, T. and Ueno, K.}, 
{\sl Monodromy preserving deformation of linear ordinary 
differential equations with rational coefficients I --- 
General theory and $\tau$-functions}, 
Physica {\bf 2D} (1981), 306--352. 
%
\bibitem{Kimura}{Kimura, H.}, 
{\sl The construction of a general solution of a Hamiltonian system 
with regular type singularity and its application to Painlev\'e 
equations}, Ann. Mat. Pura Appl. {\bf 134} (1983), 363--392. 
%
\bibitem{Malmquist}{Malmquist, J.},  
{\sl Sur les \'equations dif\'erentielles du second ordre, dont  
l'int\'egrales g\'en\'erales a ses points critiques fixes}, 
Ark. f. mat. astr. fys. {\bf 17}, No. 8 (1923), 1--89.  
%
\bibitem{Manin}{Manin, Y.I.}, 
{\sl Sixth Painlev\'e equation, universal elliptic curve, 
and mirror of $\P^2$}, 
in ``Geometry of differential equations", 
Amer. Math. Soc. Transl., Ser. 2, {\bf 186}, Amer. Math. Soc., 
Providence, 1998, 131--151.  
%
\bibitem{MMT}{Matano, T, Matumiya, A. and Takano, K.}, 
{\sl On some Hamiltonian structures of Painlev\'e systems, II},
J. Math. Soc. Japan, {\bf 51}, No.4 (1999), 843--866. 
%
\bibitem{NY}{Noumi, M. and Yamada, Y.}, 
{\sl A new Lax pair for the sixth Painlev\'e equation 
associated with} $\hat{\mathfrak{so}}(8)$, in  ``Microlocal 
analysis and complex Fourier analysis", 
T. Kawai and K. Fujita eds., World Scientific, NJ, 2002, 
238--252.   
%
\bibitem{O1}{Okamoto, K.},
{\sl Sur les feuilletages associ\'es aux \'equations du second 
ordre \`a points critiques fixes de P. Painlev\'e, Espaces des 
conditions initiales}, Japan. J. Math. {\bf 5}, (1979), 1--79. 
%
\bibitem{Okamoto1}{Okamoto, K.}, 
{\sl Polynomial Hamiltonians associated with Painlev\'e equations, I}, 
Proc. Japan Acad. {\bf56}, Ser. A (1980), 264--268; {\sl II}, ibid. 
367--371. 
%
\bibitem{Okamoto2}{Okamoto, K.}, 
{\sl Studies on the Painlev\'e equations I, sixth Painlev\'e equation 
$\PVI$}, Ann. Math. Pura Appl. (4) {\bf 146} (1987), 337--381. 
%
\bibitem{STT02} {Saito, M.-H., Takebe, T. and Terajima, H.},  
{\sl Deformationof Okamoto-Painlev\'e pairs and Painlev\'e equations},
J. Algebraic Geom. {\bf 11} (2002), no. 2, 311--362. 
%
\bibitem{SU}{Saito, M.-H. and Umemura, H.}, {\sl Painlev\'e 
equations and deformations of rational surfaces with rational 
double points}, in ``Physics and Combinatorics 1999 (Nagoya)'', 
World Sci.Publishing, River Edge, NJ, 2001, 320--365. 
%
\bibitem{STe}{Saito, M.-H. and Terajima, H.}, 
{\it Nodal curves and Riccati solutions of Painlev\'e equations},
Preprint, Kobe, 2002, (math.AG/0201225). 
%
\bibitem{Sakai}{Sakai, H.}, 
{\sl Rational surfaces associated with affine root systems and 
geometry of the Painlev\'e equations}, Comm. Math. Phys. 
{\bf 220} (2001), 165--229. 
%
\bibitem{Schlesinger}{Schlesinger, L.}, 
{\sl \"Uber eine Klasse von Differentialsystemen beliebliger 
Ordnung mit festen Kritischer Punkten}, 
J. reine angew. Math. {\bf 141} (1912), 96--145. 
%
\bibitem{Shioda-Takano}{Shioda, T. and Takano, K.}, 
{\sl On some Hamiltonian structures of Painlev\'e systems I},
Funkcial. Ekvac. {\bf 40} (1997), 271--291.
%
\bibitem{Takano}{Takano, K.}, 
{\sl Reduction for Painlev\'e equations at the fixed singular 
points of the first kind}, Funkcial. Ekvac. {\bf 29} (1986), 
99--119. 
%
\bibitem{Terajima}{Terajima, H.}, 
{\sl On the space of monodromy data of Painlev\'e VI}, 
Preprint, Kobe (2003).  
%
\bibitem{Watanabe}{Watanabe, H.}, 
{\sl Birational canonical transformations and classical solutions 
of the sixth Painlev\'e equation}, 
Ann. Scuola Norm. Sup. Pisa {\bf 27} (1998), 379--425.  
\end{thebibliography}
\end{document}